\documentclass{article}
\usepackage{blindtext}
\usepackage[T1]{fontenc}
\usepackage[utf8]{inputenc}

\usepackage{amssymb,chicago}
\usepackage{amsfonts}
\usepackage{amsmath}
\usepackage{hyperref,color}
\usepackage{amsthm}
\usepackage{geometry}
\usepackage{indentfirst}
\usepackage{bbm}

\newcommand{\eref}[1]{(\ref{#1})}
\usepackage[]{graphicx}
\usepackage{
tikz,
relsize,
booktabs,
tikz
}

\pgfdeclarelayer{background layer}
\pgfdeclarelayer{foreground layer}
\pgfsetlayers{background layer,main,foreground layer}

\tikzset{
>=latex
}
\usepackage{bbm}
\usepackage{appendix}
\usepackage{etoolbox}

\usepackage{etoc}

\newtheorem{mydefinition}{Definition}
\newtheorem{theo}{Theorem}

\newtheorem{rem}{Remark}

\newtheorem{coro}{Corollary}

\newcommand{\PXe}{P_{X, \varepsilon}}
\newcommand{\PX}{P_{X}}
\newcommand{\cA}{\mathcal{A}}
\newcommand{\cC}{\mathcal{C}}
\newcommand{\EXe}{E_{X,\varepsilon}}
\newcommand{\cG}{\mathcal{G}}
\newcommand{\EX}{E_{X}}

\newcommand{\Ee}{E_{\varepsilon}}
\newcommand{\cD}{\mathcal{D}}
\newcommand{\cB}{\mathcal{B}}
\newcommand{\cT}{\mathcal{T}}
\newcommand{\cK}{\mathcal{K}}
\newcommand{\cH}{\mathcal{H}} 
\newcommand{\cP}{\mathcal{P}} 
\newcommand{\cF}{\mathcal{F}} 
\newcommand{\cX}{\mathcal{X}} 

\DeclareMathOperator*{\argmin}{arg\,min}
\newtheorem{lemma}{Lemma}
\newenvironment{example}
{\begin{exemp}\begin{em}}
{\end{em}\end{exemp}}
\newtheorem{exemp}{Example}
\newcommand\E{\mathbb{E}}

\begin{document} 

\title{\bf Risk upper bounds for RKHS ridge group sparse estimator in the regression model with non-Gaussian and non-bounded error\thanks{Halaleh Kamari, Universit\'e Paris-Saclay, France, {\sf @}, Sylvie Huet, INRAE, France, {\sf @}, Marie-Luce Taupin, Université Evry Val d'Essonne, France, {\sf @}.}}

\author{{\sc Halaleh Kamari},
{\sc Sylvie Huet},
{\sc Marie-Luce Taupin}
}
\maketitle

\begin{abstract} 
We consider the problem of estimating a meta-model of an unknown regression model with non-Gaussian and non-bounded error.
The meta-model belongs to a reproducing kernel Hilbert space constructed as a direct sum of Hilbert spaces leading to an additive decomposition including the variables and interactions between them. The estimator of this meta-model is calculated by minimizing an empirical least-squares criterion penalized by the sum of the Hilbert norm and the empirical $L^2$-norm. In this context, the upper bounds of the empirical $L^2$ risk and the $L^2$ risk of the estimator are established. 
\end{abstract}


\section{Introduction}
Let us consider the following regression model:
\begin{equation}
\label{model}
Y=m(X)+\sigma\varepsilon,\:\sigma>0,
\end{equation} 
where the variables $X=(X_1,...,X_d)$ are independent with a known law $P_X=\bigotimes_{a=1}^d P_{a}$ on $\mathcal{X}=\prod_{a=1}^d\mathcal{X}_a$, a compact subset of $\mathbb{R}^d$. 
The number $d$ of components of $X$ may be large. The model $m$ from $\mathbb{R}^d$ to $\mathbb{R}$ maybe complex, presenting strong non-linearities, and it is assumed to be square-integrable, i.e. $m\in L^2(\mathcal{X},P_X)$.

Let $\mathcal{D}$ be the set of densities, 
\begin{align}
\label{familypi}
\mathcal{D}=\Big\{\pi_\alpha:\pi_\alpha(x)=a_\alpha\exp(-\vert x\vert^\alpha),\mbox{ with } (a_\alpha)^{-1}=\int_{\mathbb{R}}\exp(-\vert x\vert^{\alpha})dx,\:\alpha>2\Big\}.
\end{align}  
In this paper, we assume that the error term $\varepsilon$ is equal to $Z/\sigma_\alpha$, where $Z$ is a random variable with density $\pi_\alpha\in\mathcal{D}$ and $\sigma^2_\alpha$ is its variance, i.e. $\mbox{var}(Z)=\sigma_\alpha^2$.

Based on $n$ data points $\{(X_i,Y_i)\}_{i=1}^n$, a meta-model that approximates the Hoeffding decomposition of $m$ is estimated. 
This meta-model belongs to a reproducing kernel Hilbert space (RKHS), which is constructed as a direct sum of Hilbert spaces (\cite{DURRANDE201357}). 
The estimation of the meta-model is carried out via a penalized least-squares minimization allowing to select the subsets of variables $X$ that contribute to predict the output $Y$ (\cite{huet:hal-01434895}). 
 
Let us be more precise on the Hoeffding decomposition. Let $\mathcal{P}$ be the set of all the subsets of $\{1,...,d\}$ with dimension $1$ to $d$, and for all $v\in\mathcal{P}$ and $X\in\mathcal{X}$, let $X_v$ be the vector with components $X_a$ for all $a\in v$.
Let also $\vert A\vert$ be the cardinality of a set $A$ and for all $v\in\mathcal{P}$, let $m_v:\mathbb{R}^{\vert v\vert}\rightarrow \mathbb{R}$ be a function of $X_v$. Then, the Hoeffding decomposition of $m$ is written as (\cite{hoeffding1948}, \cite{Sobol1993SensitivityEF}, \cite{vaart_1998}),
\begin{align}
\label{sobol}
m(X)= m_0+\sum_{v\in\mathcal{P}}m_v(X_v),
\end{align}
where $m_0$ is a constant. 

This decomposition \eref{sobol} is unique (\cite{Sobol1993SensitivityEF}), all the functions $m_v$ are centered, and they are orthogonal with respect to $L^2(\mathcal{X},P_X)$. 

The Hoeffding decomposition of $m$ is approximated by the orthogonal projection of $m$ on a RKHS $\mathcal{H}$ which is constructed as a direct sum of Hilbert spaces (\cite{DURRANDE201357}). 

Let $\langle.,.\rangle_{\mathcal{H}}$ be the inner product in $\mathcal{H}$, and let $k$ and $k_v$ be the reproducing kernels associated with the RKHS $\mathcal{H}$ and the RKHS $\mathcal{H}_v$, respectively. The properties of the RKHS $\mathcal{H}$ insures that any function $f\in\mathcal{H}$, $f:\mathcal{X}\subset \mathbb{R}^d\rightarrow\mathbb{R}$ can be written as the following decomposition: 
\begin{align}
\label{durandhoeffintro}
f(X)=\langle f,k(X,.)\rangle_{\mathcal{H}}=f_0+\sum_{v\in\mathcal{P}}f_v(X_v),
\end{align}
where $f_0$ is a constant, and $f_v:\mathbb{R}^{\vert v\vert}\rightarrow \mathbb{R}$ is defined by, 
$$f_v(X)=\langle f,k_v(X,.)\rangle_{\mathcal{H}}.$$ 
For all $v\in\mathcal{P}$, the functions $f_v(X_v)$ are centered and for all $v, v'\in\mathcal{P}$, $v\neq v'$, the functions $f_v(X_v)$ and $f_{v'}(X_{v'})$ are orthogonal with respect to $L^2(\mathcal{X},P_X)$. Therefore, the decomposition of any function $f$ presented in Equation (\ref{durandhoeffintro}) is unique and is its Hoeffding decomposition. 

The meta-model $f^*$ that approximates the Hoeffding decomposition of $m$ is defined as follows:
\begin{equation*}
f^*=\argmin_{f \in \cH}\|m-f\|^{2}_{2} =\argmin_{f \in \cH} E_X\left(m(X)-f(X)\right)^{2}.
\end{equation*}
Since the function $f^*$ belongs to the RKHS $\mathcal{H}$, its decomposition on $\mathcal{H}$ is its Hoeffding decomposition:
\begin{equation}
\label{decfet}
f^*=f^*_0+\sum_{v\in\mathcal{P}}f^*_v.
\end{equation}
And for all $v\in\mathcal{P}$, the function $f^*_v$ in Equation \eref{decfet} approximates the function $m_v$ in Equation \eref{sobol}. 

Decomposition \eref{decfet} contains $\vert\mathcal{P}\vert$ terms $f^*_v$ to be estimated. The cardinality of $\mathcal{P}$ is equal to $2^d-1$ which may be huge since it raises very quickly by increasing $d$. In order to deal with this problem, one may estimate $f^*$ by a sparse estimator $\widehat{f}\in\mathcal{H}$. To this purpose, the estimation of $f^*$ is done on the basis of $n$ observations by minimizing an empirical least-squares criterion penalized by the sum of the Hilbert norm and the empirical norm. 
This procedure, called ridge group sparse, estimates the groups $v$ that are suitable for predicting $f^*$, and the relationship between $f^*_v$ and $X_v$ for each group $v$ (\cite{huet:hal-01434895}). The estimator so obtained is called theRKHS ridge group sparse estimator.

Several authors studied the theoretical properties of estimators similar to the RKHS ridge group sparse estimator. Let us briefly review their framework and their results.

\cite{meier2009} considered an estimator similar to the RKHS ridge group sparse estimator. Instead of adding two separate sparsity and smoothness penalties, they combine these two terms into a single sparsity and smoothness penalty. In the fixed design regression model with error $\varepsilon$ that is distributed as a sub-Gaussian random variable, they established upper bounds of the empirical risk for estimating the projection of $m$ onto the set of univariate additive functions. Afterwards, \cite{Raskutti:2012:MRS:2503308.2188398} showed (in Section 3.4. of their paper) that the convergence rate of this estimator is sub-optimal.

\cite{koltchinskii2010} considered a more general RKHS including the functions that have an additive representation over kernel spaces and obtained an estimator based on a ridge group sparse type procedure. Under a global boundedness condition, they established upper bounds on the excess risk assuming that the function $m$ has a sparse representation. A global boundedness condition means that the quantity $\sup_{f\in\mathcal{H}}\sup_{X\in\mathcal{X}}\vert f(X)\vert$ is assumed to be bounded independently of dimension $d$. 
Their results are valid for a large class of loss functions, and for distributions of the observations $Y$ such that some defined boundedness conditions on the loss functions are satisfied (see Section 2.1. of their paper). In their framework, the input variables $X$ are not assumed to be independent and there is no orthogonality assumption between the kernel spaces. Instead, the authors introduced some characteristics related to the degree of \textit{dependence} of their kernel spaces which insures \textit{almost} orthogonality between these spaces. 
Their method to derive their upper bounds relies on the elementary empirical and Rademacher process methods such as symmetrization and concentration inequalities for Rademacher processes and Bernstein type exponential bounds.

\cite{Raskutti:2012:MRS:2503308.2188398} assumed that the function $m$ has a sparse univariate additive representation, i.e. $m=\sum_{a\in S} m_a(X_a)$ for $m_a(X_a)$ being univariate functions and $\vert S\vert<d$, such that each univariate function $m_a$ lies in a RKHS $\mathcal{H}_a$. 
They used the ridge group sparse procedure to calculate the estimator of $m$, and studied the theoretical properties of  their estimator in the Gaussian regression model, i.e. $\varepsilon$ in Equation \eref{model} is distributed as a centered Gaussian random variable.
They provided upper bounds for the integrated and the empirical risks and a lower bound for the integrated risk of  their estimator over spaces of sparse additive models, including polynomials, splines and Sobolev classes.

\cite{huet:hal-01434895} studied the theoretical properties of the RKHS ridge group sparse estimator in the Gaussian regression model. They derived upper bounds with respect to the $L^2$-norm and the empirical $L^2$-norm for the distance between the true function $m$ and its estimation $\widehat{f}$ into the RKHS $\mathcal{H}$. 
 
\cite{Raskutti:2012:MRS:2503308.2188398} and \cite{huet:hal-01434895} did not assume the global boundedness condition. Instead, they assumed that each function within the unit ball of the Hilbert space $\mathcal{H}_v$ is uniformly bounded by a constant. The proof of their results 
is based on the probabilistic methods of empirical Gaussian process such as concentration inequalities and Sudakov minoration (e.g. \cite{Pisier1989}, \cite{massart2000}, \cite{geer2000empirical}, \cite{ledoux2001concentration}), as well as results on the Rademacher complexity of kernel classes (\cite{Mendelson2002}, \cite{bartlett2005}). 

In this paper, the upper bounds of the empirical $L^2$ risk and the $L^2$ risk of the RKHS ridge group sparse estimator are provided, in the regression model (see Equation \eref{model}) with non-Gaussian and non-bounded error $\varepsilon$, and by considering a quadratic loss function. 
In this case the conditions assumed in \cite{koltchinskii2010} are not satisfied, and the empirical Gaussian process methods such as concentration inequalities and Sudakov minoration can not be used.

The proof of our results requires different mathematical tools than those used in the works mentioned above:
\begin{itemize}
\item a Sudakov type minoration that is satisfied for the non-Gaussian and non-bounded random variables, 
\item a concentration bound for the lower and upper tails of a convex function of the random variables $\{\varepsilon_i\}_{i=1}^n$ that are non-Gaussian and non-bounded.
\end{itemize}  
To the best of our knowledge, in our context of regression model with non-Gaussian and non-bounded error $\varepsilon$, and with quadratic loss function, the only Sudakov type minoration which allows to obtain the same rate of convergence for the RKHS ridge group sparse estimator as in the Gaussian regression model (see \cite{huet:hal-01434895}), is the one obtained by \cite{10.2307/2374931}. The minoration obtained by \cite{10.2307/2374931} is specific to the densities $\pi_\alpha\in\mathcal{D}$ as defined in Equation \eref{familypi}. This is the reason why this class of densities is considered in this work. 

Concerning the concentration bound, it can be shown that the distribution functions associated with the densities $\pi_\alpha\in\mathcal{D}$ belong to a class of distribution functions defined by \cite{adamarticle}, for which the log-Sobolev inequality (\cite{10.2307/2373688}) is satisfied. \cite{articleshu} provided bounds for the lower and upper tails of convex functions of independent random variables which satisfy the log-Sobolev inequality. Since the distribution functions associated with the densities $\pi_\alpha\in\mathcal{D}$ satisfy the log-Sobolev inequality, the concentration inequality derived by \cite{articleshu} holds for them. 

This paper is organised as follows: The RKHS construction and the procedure for estimating a meta-model are presented in Section \ref{sec:estimator}. The theoretical properties of the RKHS ridge group sparse estimator are stated in Theorem \ref{oracle} and Corollary \ref{oracle3}. The proof of Theorem \ref{oracle} is postponed in Section \ref{prooforacle}. In Section \ref{sec:arguments} the main arguments of the proof of Theorem \ref{oracle} and motivation for the choice $\pi_\alpha$ are detailed.

\section{Meta-modelling and the RKHS ridge group sparse estimator}\label{sec:estimator}
The independency between the input variables $X$ allows to write the function $m$ according to its Hoeffding decomposition presented in Equation \eref{sobol},
$$m(X)= m_0+\sum_{v\in\mathcal{P}}m_v(X_v).$$ 
The unknown function $m$ is approximated by its orthogonal projection, denoted $f^*$, on a RKHS, denoted $\mathcal{H}$, that is constructed as a direct sum of Hilbert spaces. The RKHS $\mathcal{H}$ is associated with a so-called ANOVA kernel which is defined in order to obtain the analytical expression of the terms of the Hoeffding decomposition of the functions of $\mathcal{H}$. As $f^*$ is the orthogonal projection of $m$ on $\mathcal{H}$, each term in its decomposition is an approximation of the associated term in the Hoeffding decomposition of $m$. The construction of the RKHS $\mathcal{H}$ has been proposed by \cite{DURRANDE201357} that we recall briefly in the following.

\subsection{RKHS construction}\label{subsec:rkhs}
Let $\mathcal{X} = \mathcal{X}_{1} \times \ldots \times \mathcal{X}_{d}$ be a subset
of $\mathbb{R}^{d}$. For each $a\in \{1,\cdots,d\}$, we choose a RKHS  $\mathcal{H}_{a}$, and its associated kernel $k_{a}$ defined on the set $\mathcal{X}_{a} \subset \mathbb{R}$ such that the two following properties are satisfied:
\begin{itemize}
\item[(i)] $k_{a}:\mathcal{X}_{a} \times \mathcal{X}_{a} \rightarrow \mathbb{R}$ is $P_{a}\otimes P_{a}$ measurable,
\item[(ii)]\label{iidurr} $E_{X_{a}}\sqrt{k_{a}(X_{a}, X_{a})} < \infty$.
\end{itemize}
The property \hyperref[iidurr]{(ii)} depends on the kernel $k_a$, $a=1,...,d$ and the distribution of $X_a$, $a=1,...,d$. It is not very restrictive since it is satisfied, for example, for any bounded kernel.

The RKHS $\mathcal{H}_{a}$ can be decomposed as a sum of two orthogonal
sub-RKHS, 
$$\mathcal{H}_{a} = \mathcal{H}_{0 a}
\stackrel{\perp}{\oplus} \mathcal{H}_{1 a},$$ where $\mathcal{H}_{0 a}$ is the RKHS of zero mean functions,
\begin{align*}
 \mathcal{H}_{0 a} = \Big\{ f_{a} \in \mathcal{H}_{a},\ E_{X_{a}}(f_{a}(X_{a})) =
   0\Big\},
\end{align*}
and $\mathcal{H}_{1 a}$ is the RKHS of constant functions, 
\begin{align*}   
\mathcal{H}_{1 a} = \Big\{  f_{a} \in \mathcal{H}_{a},\: f_{a}(X_{a}) =C \Big\}.
\end{align*}
The kernel $k_{0a}$ associated with the RKHS $\mathcal{H}_{0 a}$ is defined
as follows:
\begin{equation*}
k_{0a} (X_{a},X'_{a}) = k_{a}(X_{a},X'_{a}) - 
\frac{E_{U \sim P_{a}}(k_{a}(X_{a},U))E_{U \sim P_{a}}(k_{a}(X'_{a},U))}
{E_{(U,V)\sim P_{a}\otimes P_{a}}k_{a}(U,V)}.
\end{equation*}
Let $k_{v}(X_{v}, X'_{v}) = \prod_{a \in v} k_{0a} (X_{a},X'_{a}),$ then the ANOVA kernel $k$ is defined by: 
\begin{equation*}
 k(X, X') = \prod_{a=1}^{d} 
\left(1+k_{0a}(X_{a}, X'_{a})\right) = 
1 + \sum_{v \in \mathcal{P}} k_{v}(X_{v}, X'_{v}).
\end{equation*}
For $\mathcal{H}_{v}$ being the RKHS associated with the kernel $k_{v}$, the RKHS associated with the ANOVA kernel $k$ is then defined by:
\begin{equation*}
 \mathcal{H} = \prod_{a=1}^{d}\left( \mathbbm{1} \stackrel{\perp}{\oplus}
   \mathcal{H}_{0a}\right) = \mathbbm{1} + \sum_{v \in \mathcal{P}} \mathcal{H}_{v},
\end{equation*}
where $\perp$ denotes the $L^2$ inner product.

According to this construction, any function $f \in \mathcal{H}$ satisfies the following decomposition,
\begin{align*}
f(X)=\langle f,k(X,.)\rangle_\mathcal{H}=f_0+\sum_{v\in\mathcal{P}}f_v(X_v).
\end{align*}
which is the Hoeffding decomposition of $f$.

For more background on the RKHS spaces see \cite{aronszajn50reproducing}, \cite{saitoh1988theory}, \cite{Berlinet2004ReproducingKH}.
\subsection{Approximating the Hoeffding decomposition of $m$}
\label{approxHoeffding.st}
Let $f^*\in\mathcal{H}$ be defined as follows: 
\begin{equation*}
f^*=\argmin_{f \in \cH}\|m-f\|^{2}_{2} =\argmin_{f \in \cH} E_X\left(m(X)-f(X)\right)^{2}.
\end{equation*}
The function $f^{*}=f_{0}^{*} + \sum_{v\in \cP} f_{v}^{*}$, is the approximation of $m$ on the RKHS $\mathcal{H}$, and its Hoeffding decomposition is an approximation of the Hoeffding decomposition of $m$. Therefore, according to Equation~\eref{sobol}, for all $v\in\mathcal{P}$, each function $f^{*}_{v}$ approximates the function $m_{v}$. 

The number of functions $f^{*}_{v}$ is related to the cardinality of $\mathcal{P}$, i.e. $2^d-1$, that may be huge.  The idea is to calculate a sparse estimator of $f^*$ as an estimator of $m$. To do so, the ridge group sparse procedure as proposed by \cite{huet:hal-01434895} is used that we recall in the following.
\subsection{Ridge group sparse procedure and associated estimator}\label{estimator}
Let $n$ be the number of observations. For all $v\in\mathcal{P}$, let $X_v$ be the matrix of variables corresponding to the $v$-th group, i.e. 
$$X_v=(X_{vi},i=1,...,n,v\in\mathcal{P})\in\mathbb{R}^{n\times \vert \mathcal{P}\vert}.$$ 
For any $f\in\mathcal{H}$ such that $f=f_0+\sum_{v\in\mathcal{P}}f_v$, and for some tuning parameters $\gamma_v$, $\mu_v$, $v\in\mathcal{P}$, the ridge group sparse criterion is defined as follows:
\begin{align*}
\mathcal{L}(f)=\frac{1}{n}\sum_{i=1}^n\Big(Y_i-f_0-\sum_{v\in\mathcal{P}}f_v(X_{vi}) \Big)^2
+\sum_{v\in\mathcal{P}}\gamma_v\Vert f_v\Vert_n+\sum_{v\in\mathcal{P}}\mu_v\Vert f_v\Vert_{\mathcal{H}_v},
\end{align*}
where $\Vert f_v\Vert_n$ is the empirical $L^2$-norm of $f_v$ defined by the sample $\{X_{vi}\}_{i=1}^n$ as
$$\Vert f_v\Vert_n^2=\frac{1}{n}\sum_{i=1}^nf_v^2(X_{vi}).$$
The penalty function in the criterion $\mathcal{L}(f)$ is the sum of the Hilbert norm and the empirical norm, which allows to select few terms in the additive decomposition of $f$ over sets $v \in \mathcal{P}$. Moreover, the Hilbert norm favours the smoothness of the estimated $f_v$, $v \in \mathcal{P}$. 

Let us define the set of functions,
\begin{equation}
 \cF = \Big\{ f: f=f_{0} + \sum_{v\in \cP} f_{v},\mbox{ with }
   f_v \in \cH_{v},\mbox{ and } \|f_v\|_{\cH_v} \leq r_v,\: r_v>0 \Big\}.
\label{calF.eq}
\end{equation}
Then the RKHS ridge group sparse estimator is defined by,
\begin{equation}
\label{prediction}
\widehat{f}=\argmin_{f\in\mathcal{F}}\mathcal{L}(f).
\end{equation}

\section{Risk upper bounds}
In this Section, the upper bounds of the empirical $L^2$ risk and the $L^2$ risk of the RKHS ridge group sparse estimator are presented in Theorem \ref{oracle} and Corollary \ref{oracle2}, respectively.
Before stating these results, let us introduce some notation and assumptions that are needed in the rest of this paper. 

For a function $f\in\mathcal{H}$, let $S_{f}$ be its support, 
\begin{equation}
S_{f}  = \{ v \in \mathcal{P}:\: f_v \neq 0 \}.
\label{sparsity}
\end{equation}
The RKHS construction as described in Section \ref{subsec:rkhs} insures that the following properties are satisfied: 
\begin{itemize}
\item for all $v\in\mathcal{P}$, the functions $f_v \in \mathcal{H}_v$ are centered and are square-integrable, i.e.
$$E_X (f_{v}(X_v))= 0 \mbox{ and }E_X (f^{2}_{v}(X_v))< \infty,$$  
\item for all $v, v'\in\mathcal{P}$, $v\neq v'$, the functions $f_v \in \mathcal{H}_v$ and $f_{v'} \in \mathcal{H}_{v'}$ are orthogonal with respect to $L^2(\mathcal{X},P_X)$, i.e.
 $$  E_X (f_{v}(X_v) f_{v'}(X_{v'})) = 0.$$
\end{itemize}
We assume moreover that, 
\begin{itemize}
\item for all $v\in\mathcal{P}$, the functions $f_v \in \mathcal{H}_v$ are uniformly bounded, i.e. 
\begin{equation*}
\exists R>0\mbox{ such that }\Vert f_v\Vert _\infty =\sup_{X_v}\vert f_v(X_v)\vert\leq R.
\end{equation*}
\end{itemize}
Each kernel $k_v$, $v \in \mathcal{P}$ is associated with an integral operator $T_{k_v}$ from $L^2(\mathcal{X}_v,P_{v})$ to $L^2(\mathcal{X}_v,P_{v})$ defined by: 
$$\forall f\in L^2(\mathcal{X}_v,P_{v}),\: T_{k_v}(f)=\int_{\mathcal{X}_v}k_v(.,t)f(t)dP_v(t).$$
For each $v \in \mathcal{P}$, let  $\omega_{v,1}\geq \omega_{v,2}\geq ...\geq 0$ be the eigenvalues of the integral operator $T_{k_v}$ (see Equation \eref{kerneloperator}). Let us define the function $Q_{n,v}(t)$ for some positive $t$ as follows:
\begin{equation}
\label{Qn}
Q_{n,v}(t)=\sqrt{\frac{5}{n} \sum_{\ell\geq 1} \min(t^2,  \omega_{v,\ell}  )},
\end{equation}
and for some $\Delta > 0$  let $\nu_{n,v}$ be defined by:
\begin{align}
\label{nu}
\nu_{n, v} = \inf_t\Big\{Q_{n,v}(t)\leq  \Delta t^{2}\Big\}.
\end{align}
For each $v\in \mathcal{P}$, $\nu_{n,v}$ refers to the minimax optimal rate for $L^2(\mathcal{X},P_X)$-estimation in the RKHS $\mathcal{H}_v$ (\cite{Mendelson2002}). 
\begin{rem}
The rate $\nu_{n,v}$, $v\in\mathcal{P}$, depends on the regularity of the RKHS via the decreasing rate of the eigenvalues $\{\omega_{v,\ell}\}_{\ell=1}^\infty$. When RKHS is of high regularity, i.e. when the eigenvalues $\{\omega_{v,\ell}\}_{\ell=1}^\infty$ decrease quickly, then the rate $\nu_{n,v}$, $v\in\mathcal{P}$ will be close to the parametric rate of convergence (see Section \ref{RateCvge.st}).  
\end{rem}
The choice of tuning parameters in the criterion $\mathcal{L}(f)$ is specified in terms of the following quantity:
\begin{align}
 \lambda_{n,v} = \max\Big(
  \nu_{n,v}, \sqrt{\frac{d}{n}} \Big). \label{lambda.eq}
\end{align}

\begin{theo}
\label{oracle}
Consider the regression model defined at Equation~\eref{model} with $\sigma=1$.
Let  $\{(Y_{i}, X_{i})\}_{i=1}^n$ be a $n$-sample with the same law as $(Y,X)$, and let $\{\varepsilon_i\}_{i=1}^n$ be the random errors that are independent and identically distributed (i.i.d.) like $\varepsilon$.
Let also $\widehat{f}$ be defined by~\eref{prediction} with $r_v=1$ in \eref{calF.eq}, and let the tuning parameters $\mu_{v}$'s and $\gamma_{v}$'s be chosen as follows: 

For some constant $C_{1}>10+4\Delta$, 
\begin{equation}
\label{condmu}
\forall v \in \cP ,\: \mu_{v} = C_{1}\lambda^{2}_{n,v} ,\: \gamma_{v} = C_{1} \lambda_{n,v}.
\end{equation}

If there exists positive constants $C_{2}, C_{3}$, and $0<\beta<1/\alpha$ such that
 the following assumptions are satisfied:
\begin{equation}
\forall v \in \cP, \: 
n\lambda_{n,v}^{2} \geq -C_{2} \log\lambda_{n,v},
\end{equation}
and 
\begin{equation}
\label{hypdf.eq}
\forall f \in \cF,  \: 
\sum_{v \in S_{f}} \lambda_{n,v}^{2} \leq C_{3} n^{2\beta -1},
\end{equation}
then, there exists $0 < \eta < 1 $ depending on constants $\{C_i\}_{i=1}^3$, $\beta$, and $n$ ($\eta$ tends to $0$ as $n$ increases), such that with probability greater than $1-\eta$, we have for some constant $C$,
\begin{equation}
\label{ch1riskbound}
 \| m - \widehat{f}\|^{2}_{n} \leq
C \inf_{f \in \cF}\Big\{  \| m - f\|^{2}_{n}
+  \sum_{v \in S_{f}} (\mu_{v} + \gamma_{v}^{2})
 \Big\}.
\end{equation}
\end{theo}
Let us now comment on the theorem.
\begin{rem}
Let $f'$ be the function in $\mathcal{F}$ such that the infimum of the right hand side of the inequality \eref{ch1riskbound} is realized. The term $\|m-f'\|_n^2$ is the usual bias term. It quantifies 
both the approximation properties of the RKHS $\mathcal{H}$, and the bias-variance 
trade-off.
\end{rem}
\begin{rem}
This result is similar  to the  one obtained in the Gaussian regression model at the cost of the additional Assumption \eref{hypdf.eq}. This assumption allows to obtain the same rate of convergence for the RKHS ridge group sparse estimator as in the Gaussian regression model (see \cite{huet:hal-01434895}). However, it implies some restrictions on the regularity of the RKHS $\mathcal{H}$. Indeed, as for all $v\in\mathcal{P}$, $\lambda_{n,v}\geq\nu_{n,v}$ (see Equation \eref{lambda.eq}), it follows that $\sum_{v\in S_f} \nu_{n,v}^{2} \leq C_3 n^{2\beta-1}$, which implies some restrictions on the regularity of the RKHS: if $\beta$ is small, which will be the case if $\alpha$ is
large, then the RKHS should be of high regularity.
\end{rem}
\begin{rem}\label{secondlambda} 
By Equation \eref{lambda.eq}, we also have that for all $v\in\mathcal{P}$, $\lambda_{n,v}\geq \sqrt{d/n}.$ This assumption allows to control the probability of the $\vert\mathcal{P}\vert$ events (see Equation \eref{conrtolproba}), where $\log(\vert\mathcal{P}\vert)$ is of order $d$. 
\end{rem}
\begin{rem}\label{ch2:remarkfr}
The result in Theorem \ref{oracle} can be generalized to the case where $\sigma\neq1$ in Equation \eref{model}, and where $r_v\neq1$ in \eref{calF.eq}.  

Let $\widehat{g}$ be defined as follows: 
\begin{equation}
\label{ch2:ghat}
\widehat{g} = \argmin_{g\in\mathcal{F}'}\Big\{ 
\| \frac{Y}{\sigma}  - g\|_{n}^{2} +
\frac{1}{\sigma}\sum_{v}\gamma_{v}\|g_{v}\|_{n}  + 
\frac{1}{\sigma}\sum_{v}\mu_{v}\|g_{v}\|_{\mathcal{H}_{v}}
\Big\},
\end{equation}
with 
\begin{equation}
\label{ch2:fprim}
\mathcal{F}'=\Big\{g:g=g_0+\sum_{v} g_{v},\mbox{ with }g_v  \in \mathcal{H}_{v},\mbox{ and } \|g_{v}\|_{\mathcal{H}_{v}} \leq \frac{r_v}{\sigma}\Big\}.
\end{equation}
We have $\widehat{f}= \sigma \widehat{g}$ for $\widehat{f}$ being defined by \eref{prediction}.

\medskip

For all $u>0$, let $\mathcal{H}_{v}^{u}$ be the RKHS associated with the kernel $u k_{v}$. If  $u=r_{v}^{2}/\sigma^{2}$, then
 \begin{equation*}
\widehat{g} = \argmin_{g\in\mathcal{F}''}\Big\{ 
\| \frac{Y}{\sigma}  - g\|_{n}^{2} +
\frac{1}{\sigma}\sum_{v}\gamma_{v}\|g_{v}\|_{n}  + 
\frac{1}{\sigma^{2}}\sum_{v}\mu_{v} r_{v}
\|g_{v}\|_{\mathcal{H}^{u}_{v}}\Big\}. 
\end{equation*}
where
\begin{equation*}
\mathcal{F}''=\Big\{g:g=g_0+\sum_{v} g_{v},\mbox{ with }g_v  \in \mathcal{H}^{u}_{v},\mbox{ and }\|g_{v}\|_{\mathcal{H}^{u}_{v}} \leq 1\Big\}.
\end{equation*}

\medskip

We apply  Theorem \ref{oracle} with $Y/\sigma$ and $m/\sigma$ in place of $Y$
and $m$, to $\widehat{g}$ defined as above.

Let 
\begin{equation*}
 Q^{u}_{n,v}(t) = 
\sqrt{\frac{5}{n} \sum_{\ell\geq 1} \min(t^{2},u \omega_{v,\ell})},
\end{equation*} and for $\Delta' >0$, let
\begin{equation*}
\nu^{u}_{n,v}(\Delta') = \inf_t\Big\{ Q^{u}_{n,v}(t) \leq
  \Delta' t^{2}\Big\}.
\end{equation*}
Let also
\begin{equation*}
\lambda^{u}_{n,v}  = \max \Big(\nu^{u}_{n,v} , \sqrt{\frac{d}{n}}\Big).
\end{equation*}
For some constant $C_{1} > 10 + \Delta'$, take 
\begin{equation*}
 \frac{\mu_{v} r_{v}}{\sigma^{2}} = C_{1}\Big(
   \lambda^{u}_{n,v}\Big)^{2}, \: 
\frac{\gamma_{v}}{\sigma} = C_{1} \lambda^{u}_{n,v}.
\end{equation*}
Then, for $S_g$ being defined as follows
\begin{equation}
\label{ch2:sg}
S_g=\{v\in\mathcal{P}:g_v\neq0\},
\end{equation}
we have
\begin{equation*}
 \| \frac{m}{\sigma}  -  \widehat{g}\|_{n}^{2}
\leq C \inf_{g\in\mathcal{F}''} \Big\{ \|\frac{m}{\sigma}  -
g\|_{n}^{2}
+\frac{1}{\sigma^{2}}\sum_{v\in S_{g}} ( \mu_{v}r_{v} + \gamma_{v}^{2} )
\Big\},
\end{equation*}
or, multiplying both sides by $\sigma^{2}$, and taking $u=r_{v}^{2}/\sigma^{2}$,
\begin{equation*}
  \| m  -  \sigma \widehat{g}\|_{n}^{2}
\leq C \inf_{g\in\mathcal{F}'} \Big\{ \|m  -
\sigma g\|_{n}^{2}
+\sum_{v\in S_{g}} ( \mu_{v}r_{v} + \gamma_{v}^{2} )
\Big\}.
\end{equation*}
\end{rem}
\begin{coro}
\label{oracle2}
Under the same assumptions as Theorem \ref{oracle}, we have with high probability for some constant $C'$ that,
\begin{align*}
 \| m - \widehat{f}\|^{2}_{2} \leq
C' \inf_{f \in \cF}\Big\{  \| m - f\|^{2}_{n} + \| m - f\|^{2}_{2}
+  \sum_{v \in S_{f}} (\mu_{v} + \gamma_{v}^{2})
 \Big\}.
\end{align*}
\end{coro}
\begin{rem}\label{ch2:remarkfr2}
The result in Corollary \ref{oracle2} can be generalized to the case where $\sigma\neq1$ in Equation \eref{model}, and where $r_v\neq1$ in \eref{calF.eq}. It suffices to apply Corollary \ref{oracle2} with $Y/\sigma$ and $m/\sigma$ in place of $Y$ and $m$, to $\widehat{g}$ as defined in Equation \eref{ch2:ghat}. Then, with similar demonstration as in Remark \ref{ch2:remarkfr} we obtain,
\begin{align*}
 \| m - \sigma\widehat{g}\|^{2}_{2} \leq
C' \inf_{g \in \mathcal{F}'}\Big\{  \| m - \sigma g\|^{2}_{n} + \| m - \sigma g\|^{2}_{2}
+  \sum_{v \in S_{g}} (\mu_{v}r_v + \gamma_{v}^{2})
 \Big\},
\end{align*}
where $\mathcal{F}'$ and $S_g$ are defined in Equations \eref{ch2:fprim} and \eref{ch2:sg}, respectively. 
\end{rem}  
\subsection{Rate of convergence}
\label{RateCvge.st}
\begin{coro}
\label{oracle3}
Under the same assumptions as Theorem~\ref{oracle}, we have
\begin{equation*}
 \| m - \widehat{f}\|^{2}_{n} \leq
C \inf_{f \in \cF}\Big\{  \| m - f\|^{2}_{n}
+ \Big(\sum_{v \in S_{f}} \nu^{2}_{n,v} + \frac{d |S_{f}|}{n}\Big) \Big\}.
\end{equation*}
\end{coro} 
This Corollary highlights that the upper bound is relevant when the infimum is reached for functions $f$ that have a sparse decomposition in $\cH$, i.e. $\vert S_f\vert$ is small, and when $d$ is small face to $n$. When $d$ is large, the  decomposition of functions in $\cH$ should be limited to interactions of a limited order, so that the number of 
elements in the estimated meta-model is of order smaller than $d^r$ for some small $r$, say $r=2$ for example. 
In such a case, the cardinality of $\mathcal{P}$ will be smaller than $d^2$. 
As we mentioned in Remark \ref{secondlambda}, the assumption $\lambda_{n,v}\geq \sqrt{d/n}$ is needed to control the value $\log(\vert \mathcal{P}\vert)$, which will be now smaller than $2\log(d)$. Therefore, the value $d$ in the definition of $\lambda_{n,v}$ (see Equation \eref{lambda.eq}) as well as the term $d\vert S_f\vert/n$ in the infimum above will be replaced by $2\log(d)$ and $2\log(d)\vert S_f\vert/n$, respectively.

Let us  discuss the rate of convergence given by
$\sum_{v \in S_{f}} \nu^{2}_{n, v}$. 
For the sake of simplicity we consider the case where the
variables $X_{1}, \ldots, X_{d}$ have the same distribution $P_{1}$
on $\cX_{1} \subset\mathbb{R}$, and
where the 
unidimensional kernels $k_{0 a}$ are all identical, such that
$k_{v}(X_{v}, X'_{v}) = \prod_{a \in v}
k_{0}(X_{a},X'_{a})$. The kernel $k_{0}$ admits an eigen expansion
given by 
\begin{equation*}
 k_{0}(X_a,X'_a) = \sum_{\ell_a \geq 1} \omega_{0, \ell_a} \phi_{\ell_a}(X_a)\phi_{\ell_a}(X'_a), 
\end{equation*}
where the eigenvalues $\{\omega_{0, \ell_a}\}_{\ell_a=1}^\infty$ are non-negative and ranged
in the decreasing order, and where the $\{\phi_{\ell_a}\}_{\ell_a=1}^\infty$ are the
associated eigenfunctions, orthonormal with respect to $L^{2}(\mathcal{X}_1, P_1)$.
Therefore, the kernel $k_{v}$  admits the following expansion,
\begin{equation}
\label{kerneloperator}
 k_{v}(X_v,X^{'}_{v}) = \sum_{\ell=(\ell_1 \ldots \ell_{|v|})} 
\underbrace{\prod_{a=1}^{|v|}\omega_{0, \ell_{a}}}_{\omega_{v, \ell}}
\underbrace{\prod_{a=1}^{|v|}\phi_{\ell_{a}}(X_{a})}_{\phi_{v, \ell}(X_{v})}
\underbrace{\prod_{a=1}^{|v|}\phi_{\ell_{a}}(X^{'}_{a})}_{\phi_{v, \ell}(X^{'}_{v})}.
\end{equation}
Consider the case where the eigenvalues $\{\omega_{0, \ell_a}\}_{\ell_a=1}^\infty$ are
decreasing at a rate $\ell_a^{-2\alpha'}$ for some $\alpha' > 1/2$, i.e. the $\omega_{0, \ell}$ are of order $\ell^{-2\alpha'}=(\prod_{a=1}^{|v|}\ell_a)^{-2\alpha'}$. 
It is shown in Section 8.3. of \cite{huet:hal-01434895}, that 
$$\nu_{n, v}\propto n^{-\frac{\alpha'}{(2\alpha' + 1)}} (\log n)^{\gamma'},$$ where the rate $\nu_{n, v}$ is defined at Equation~\eref{nu} and $$\gamma' \geq
(|v|-1)\frac{\alpha'}{(2\alpha'-1)}.$$ 
For all $f \in \cF$ we have then,
\begin{equation*}
\sum_{v \in S_{f}}\nu_{n, v}^{2} \propto |S_{f}|
n^{-\frac{2\alpha'}{(2\alpha' + 1)}}(\log n)^{2\gamma'}.
\end{equation*}
Note that in this particular case, the rate of convergence depends on $|v|$ through the logarithmic term $(\log n)^{2\gamma'}$, and that up to this logarithmic term the rate of convergence has the same order than the usual non-parametric  rate for unidimensional functions. It follows that the RKHS space $\cH$ should be chosen such that the unknown function $m$ is well approximated by sparse functions in $\cH$ with low order of interactions.

Besides, the rate $\nu_{n, v}$ should satisfy assumption \eref{hypdf.eq}, 
\begin{equation*}
 \sum_{v \in S_{f}}\nu_{n, v}^{2} \leq C_{3} n^{2\beta-1},
\end{equation*}
which holds if
\begin{equation}
\label{alpp}
\alpha'  > \frac{1-2\beta}{4\beta} > \frac{\alpha-2}{4}.
\end{equation}
This shows that for the large values of $\alpha$ the assumption \eref{hypdf.eq} implies some restrictions on the regularity of the RKHS chosen:  If $\alpha<4$, then all $\alpha'$ greater than $1/2$ satisfy Equation \eref{alpp}, since $(\alpha-2)/4 < 1/2.$ If $\alpha\geq4$, then we have $\alpha' >(\alpha-2)/4>1/2$. As $\alpha$ increases, i.e. $\beta$ decreases (recall that $0<\beta<1/\alpha$), and assumption \eref{hypdf.eq} implies that the RKHS chosen should be of high regularity.


\section{Main arguments of the proof of Theorem \ref{oracle} and motivation for the choice $\pi_\alpha$}\label{sec:arguments}
The proof of Theorem \ref{oracle} starts in the same way as the proof
of Theorem $2.1.$ in \cite{huet:hal-01434895} where they considered the Gaussian regression model. However, it differs in two essential points:
\begin{itemize}
\item[1.] Sudakov type minoration,
\item[2.] Concentration inequality.
\end{itemize}
In the following Section, we give a sketch of the proof of Theorem~\ref{oracle}, we highlight the two points above that differs the proof from the proof in the Gaussian regression model, and we provide a detailed comparison to the related works. 
In Section \ref{sec:sudakov} we give a brief introduction to the Sudakov type minoration context, we explain the motivation for choosing densities $\pi_\alpha\in\mathcal{D}$ defined in Equation \eref{familypi}, and we state in Corollary \ref{corosudakov} the appropriate Sudakov minoration used in the proof of Theorem \ref{oracle}. In Section \ref{sec:conineq} we present the concentration inequality context, and we state in Corollary \ref{lipschitz} the appropriate concentration inequality used in the proof of Theorem \ref{oracle}.  

\subsection{Sketch of the proof} \label{sketch}
We give here a sketch of the proof of Theorem \ref{oracle}, and we postpone to Section \ref{prooforacle} for complete statements. We begin by introducing some notation. 

We denote by $C$ constants that vary from an equation to the other.
For $v\in\mathcal{P}$, and for a function $\phi :\mathbb{R}^{|v|}\mapsto \mathbb{R}$, we denote by $V_{n,\varepsilon}$ the empirical process defined as,
\begin{align}
\label{procemp}
V_{n,\varepsilon}(\phi)=\frac{1}{n}\sum_{i=1}^n \varepsilon_i \phi(X_{v,i}).
\end{align} 
For all $v\in \mathcal{P}$, let $\cH_{v}$ be the RKHS associated with the reproducing kernel $k_v$. For any function $g_v \in \cH_{v}$, $v\in\mathcal{P}$, and $V_{n, \varepsilon}$ being defined in Equation~\eref{procemp}, we consider two following processes,
\begin{align}
 W_{n,2,v}( t) &= \sup \Big\{|V_{n, \varepsilon}(g_{v})|, \:
\|g_v\|_{\mathcal{H}_v}\leq 2,\: \|g_v\|_{2}\leq t\Big\} \label{defWn}, \\
W_{n,n,v}( t) &= \sup \Big\{ |V_{n, \varepsilon}(g_{v})|,\: 
\|g_v\|_{\mathcal{H}_v}\leq 2, \: \|g_v\|_{n}\leq t\Big\}.\label{defWntilde}
\end{align} 
Starting from the definition of $\widehat{f}$, some simple calculations give that for all $f\in \cF$, 
\begin{align*}
 C \|m - \widehat{f}\|_{n}^{2} \leq &\|m - f\|_{n}^{2} +|V_{n,\varepsilon}(\widehat{f}-f)|+ \sum_{v \in S_{f}}[\gamma_{v} \|\widehat{f}_{v}-f_{v}\|_{n} + \mu_{v}\|\widehat{f}_{v} - f_{v}\|_{\cH_{v}}]\\
 &-\sum_{v\notin S_f}[\mu_v\| \widehat{ f}_v\|_{\mathcal{H}_v}+  \gamma_v\| \widehat{f}_v\|_{n}],\\
\leq &\|m - f\|_{n}^{2} +|V_{n,\varepsilon}(\widehat{f}-f)|+ \sum_{v \in S_{f}}[\gamma_{v} \|\widehat{f}_{v}-f_{v}\|_{n} + \mu_{v}\|\widehat{f}_{v} - f_{v}\|_{\cH_{v}}].
\end{align*} 
If we set $g=\widehat{f} - f$, then $g \in\cH$, $g=g_{0} + \sum_{v} g_{v}$, with $g_{v} = \widehat{f}_v - f_v$, and for each $v$, $\|g_{v}\|_{\cH_{v}}
\leq 2$.

The main problem is now to control the empirical process $V_{n,\varepsilon}$.  
For each $v$, letting $\lambda_{n,v}$ as in (\ref{lambda.eq}),
we state (see Lemma~\ref{Tau}, page~\pageref{Tau}) that, with high
probability, 
\begin{align}
\label{lem1sketch}
\vert V_{n,\varepsilon}(g_v)\vert\leq C\lambda_{n,v}^2\Vert g_v\Vert_{\mathcal{H}_v}+C\lambda_{n,v}\Vert g_v\Vert_n.
\end{align}

One of the key points in the proof of Lemma \ref{Tau} is to find an upper bound for the two following quantities:
\begin{align}
\label{cncntsketch}
\vert W_{n,n,v}(t)-E_\varepsilon(W_{n,n,v}(t))\vert,\mbox{ and }\vert W_{n,2,v}(t)-E_\varepsilon(W_{n,2,v}(t))\vert .
\end{align}
In the Gaussian regression model, one use the isoperimetric inequality for Gaussian processes in \cite{massart2007concentration}. 

When dealing with errors that are not distributed as a Gaussian distribution, different tools are needed to obtain the upper bounds for the quantities in Equation \eref{cncntsketch} (see Section \ref{sec:conineq} for a complete discussion of this point of the proof). 
Let us continue the sketch of the proof before coming back to this point. 

If for all $v$, $\mu_{v}$ and  $\gamma_{v}$ satisfying Equation~\eref{condmu}, by using Equation \eref{lem1sketch} we deduce that with high probability,
\begin{equation*}
C \|m - \widehat{f}\|_{n}^{2} \leq   \|m - f\|_{n}^{2}+
\sum_{v \in S_f} [\gamma_{v}\|g_{v}\|_{n} + \mu_{v}
\|g_{v}\|_{\cH_{v}}]
+ \sum_{v \notin S_{f}}
[\gamma_{v} \|\widehat{f}_{v}\|_{n} + \mu_{v} \|\widehat{f}_{v}\|_{\cH_{v}}].
\end{equation*}
Besides, we can express the decomposability property of the penalty as
follows (see lemma~\ref{lemme2}, page~\pageref{lemme2}): 

over the set where the empirical process  is controlled as stated above, we have with high probability, 
\begin{equation*}
\sum_{v \notin S_{f}}[\gamma_{v} \|\widehat{f}_{v}\|_{n} + \mu_{v}\|\widehat{f}_{v}\|_{\cH_{v}}]
\leq  C \sum_{v \in S_{f}}[\gamma_{v} \|g_{v}\|_{n} + \mu_{v}\|g_{v}\|_{\cH_{v}}].
\end{equation*}
Putting the things together, and using that $\|g_{v}\|_{\cH_{v}} \leq
2$, we obtain the following upper bound:
\begin{equation*}
 C \|m - \widehat{f}\|_{n}^{2} \leq  \|m - f\|_{n}^{2} + \sum_{v \in S_{f}} [\mu_{v} +
\gamma_{v}  \|g_{v}\|_{n}].
\end{equation*}
The last important step consists in comparing $\sum_{v \in S_{f}} \|g_{v}\|_{n} $ to
$\|\sum_{v \in S_{f}} g_{v}\|_{n}$. To do so, we show first (see lemma~\ref{norm2normn} page~\pageref{norm2normn}) that for all $v\in \cP$, with high probability,
\begin{equation*}
\|g_{v}\|_{n}  \leq 2 \|g_{v}\|_2 + \gamma_{v}.
\end{equation*}
Using inequality above and that for all positive $K$, $2ab\leq (1/K)a^2+Kb^2$ we obtain,
\begin{align*}
 C \|m - \widehat{f}\|_{n}^{2} &\leq  \|m - f\|_{n}^{2}+ \sum_{v \in S_{f}} (\mu_{v}
  + \gamma^{2}_{v}) + \sum_{v\in S_f}\Vert g_v\Vert_{2}^2,\\
&\leq  \|m - f\|_{n}^{2}+ \sum_{v \in S_f} (\mu_{v}
  + \gamma^{2}_{v}) + \sum_{v\in \mathcal{P}}\Vert g_v\Vert_{2}^2.  
\end{align*}
Then we use the orthogonality assumption between the spaces $\cH_{v}$,
$$\sum_{v\in \mathcal{P}}\Vert g_v \Vert^2_{2}=\Vert \sum_{v\in \mathcal{P}}g_{v}\Vert^2_{2}= \Vert g\Vert_2^2,$$ 
which allows us to obtain the following result:
\begin{equation*}
 C \|m - \widehat{f}\|_{n}^{2} \leq  \|m - f\|_{n}^{2}+ \sum_{v \in S_{f}} (\mu_{v}
  + \gamma^{2}_{v}) + 
\Vert\widehat{f} - f\Vert^{2}_{2}.
\end{equation*}
It remains now to consider different cases according to the
rankings of $\|\widehat{f}-f\|^{2}_{2}$ and
$\|\widehat{f}-f\|^{2}_n$  to get the result of Theorem~\ref{oracle}. 

If $\Vert\widehat{f}-f\Vert_2\leq\Vert\widehat{f}-f\Vert_n$ the result is obtained by a simple rearrangement of the terms. 

If $\Vert\widehat{f}-f\Vert_2\geq\Vert\widehat{f}-f\Vert_n$, under some suitable assumptions it is shown (see Lemma \ref{lemcompnormes3} page \pageref{lemcompnormes3}) that with high probability we have 
$$\Vert\widehat{f}-f\Vert_2\leq\sqrt{2}\Vert \widehat{f}-f\Vert_n.$$ 
One of the steps to prove the inequality above is to lower bound the expectation of the supremum of the empirical process, i.e. $\mathbb{E}_\varepsilon \sup_g \vert V_{n,\varepsilon}(g)\vert$ by a function of the covering number of the functional class under study, say $\mathcal{G}$.
In order to solve this step in the Gaussian regression model one may use the Sudakov minoration in \cite{Pisier1989}, for which the minoration is obtained thanks to the Slepian’s Lemma. 
The Slepian’s Lemma is specific to the Gaussian setting, and it does not hold when dealing with errors that are not distributed as a centered Gaussian distribution.

In the regression model (see Equation \eref{model}) with error $\varepsilon$ that is distributed with density proportional to $\pi_\alpha\in\mathcal{D}$, the proof of the upper bound stated in Theorem \ref{oracle}, needs two following mathematical tools:
\begin{itemize}
\item[Point 1.] a Sudakov type minoration to link the covering number on a class $\mathcal{G}$ to the expectation of the supremum of the empirical process over this class $\mathcal{G}$, $\mathbb{E}_\varepsilon \sup_{g\in\mathcal{G}} \vert V_{n,\varepsilon}(g)\vert$, and conclude Lemma \ref{lemcompnormes3},
\item[Point 2.] a concentration inequality to bound the quantities defined in Equation \eref{cncntsketch} which leads to bound the empirical process $V_{n,\varepsilon}$ and conclude Lemma \ref{Tau}.
\end{itemize}
The Point $1.$ is solved using a Sudakov type minoration which is a consequence of the result obtained by \cite{10.2307/2374931}. More precisely, it can be shown (see Corollary \ref{corosudakov} page \pageref{corosudakov}) that for $\varepsilon=(\varepsilon_1,...,\varepsilon_n)$ being i.i.d. random variables distributed with density $\pi_\alpha\in\mathcal{D}$ (see Equation \eref{familypi}), and for all $\delta>0$, we have,
\begin{align}
\label{sudtalagrand}
\frac{1}{K}\log N(\delta,\mathcal{G},\Vert.\Vert) &\leq (\frac{2nE_{\varepsilon}\sup_{g\in\mathcal{G}}\vert V_{n,\varepsilon}(g)\vert}{\delta})^2 \mathbf{1}_{[2nE_{\varepsilon}\sup_{g\in\mathcal{G}}\vert V_{n,\varepsilon}(g)\vert,\infty)}(\delta)\nonumber\\ 
&+(\frac{2nE_{\varepsilon}\sup_{g\in\mathcal{G}}\vert V_{n,\varepsilon}(g)\vert}{\delta})^{\alpha}\mathbf{1}_{(0, 2nE_{\varepsilon}\sup_{g\in\mathcal{G}}\vert V_{n,\varepsilon}(g)\vert]}(\delta),
\end{align}
where $K$ is a constant that depends on $\alpha$ only, $\Vert.\Vert$ is the Euclidean norm, $N(\delta,\mathcal{G},\Vert.\Vert)$ is the $\delta$-covering number 
of the metric space $(\mathcal{G},\Vert.\Vert)$, and $1_A:\mathcal{A}\rightarrow\{0,1\}$ is the indicator function of $A\subset\mathcal{A}$, i.e.
\begin{align*}
1_A(a)=\left\{ \begin{array}{rcl}
         1 & \mbox{if}
         & a\in A, \\ 
         0  & \mbox{if} & a\notin A .
                \end{array}\right.
\end{align*}

The proof of Lemma \ref{lemcompnormes3} proceeds using Equation \eref{lem1sketch} and is concluded under the Hypothesis \eref{condmu} and \eref{hypdf.eq}.

The Point $2.$ is solved using a concentration inequality (see Corollary \ref{lipschitz} page \pageref{lipschitz}) which is a consequence of the result obtained by \cite{articleshu}. 

\hfill $\Box$ 

The appropriate results to solve Point $1.$ and Point $2.$ are stated in Corollary \ref{corosudakov} in Section \ref{ssudakov} and Corollary \ref{lipschitz} in Section \ref{sconineq}, respectively. 
\subsubsection{Comparison with related works}\label{comparison}
\begin{itemize}
\item \cite{meier2009} considered a least-squares criterion penalized by a penalty function similar to the one we consider in our work. Their estimator of the unknown function $m$ has an univariate additive decomposition, i.e. decomposition \eref{sobol} limited to the main effects. 

They used a \textit{compatibility condition} to compare the sum of the empirical $L^2$-norm of the univariate functions to the empirical $L^2$-norm of the sum of the univariate functions. More precisely,

Let $S^*=\{a\in\{1,...,d\}, \Vert f_a\Vert_n\neq0\}$, then for $C(f_a)$ being a term depending on the functions $f_a$, $a\in S^*$,
\begin{align*}
\sum_{a\in S^*} \|f_a\|_{n}^2\leq \|\sum_{a=1}^df_a \|_{n}^2+C(f_a).
\end{align*}

The control of the Empirical process is done in their Lemma 1. This lemma is proved using Lemma 8.4 in \cite{geer2000empirical}, for which the errors should have sub-Gaussian tails, i.e.
\begin{align*}
\max_i E\Big(\exp(\frac{\varepsilon_i^2}{C_1})\Big)\leq C_2,
\end{align*} 
where $C_1$ and $C_2$ are constants.

Afterwards, it was shown by \cite{Raskutti:2012:MRS:2503308.2188398} (see Section 3.4. of their paper) that the convergence rate of this estimator is sub-optimal.
\item \cite{koltchinskii2010} 
considered a large class of loss functions, called \textit{loss functions of quadratic type}, which satisfies the boundedness conditions. More precisely, for $l$ being a loss function, they assume that $l(Y,.)$ is uniformly bounded from above by a numerical constant. So for a given distribution of the observations $Y$, there may exists a loss function that belongs to the class of the \textit{loss functions of quadratic type} (see Section 2.1. of their paper for some examples). 

They consider the input variables $X$ that may be not independent, and they do not assume that there is orthogonality between their RKHS, therefore $\Vert \sum_v f_v\Vert_2\neq\sum_v\Vert  f_v\Vert_2$. Instead, in their Section 2.2., they introduce some geometric characteristics related to the degree of \textit{dependence} of their RKHS, which insures \textit{almost} orthogonality between these spaces.
  
The control of the empirical process is done in their Lemma 9. This lemma is proved under the global boundedness condition and the assumptions of the \textit{loss functions of quadratic type}.  
 
We consider the quadratic loss function to obtain an estimator of the function $m$ in the regression model defined in Equation \eref{model}, with error $\varepsilon$ that is non-bounded. This case is not included in the class of the \textit{loss functions of quadratic type}. We do not impose the global boundedness condition. Instead, we assume that for all $v\in\mathcal{P}$ the functions $f_v$ are uniformly bounded. More precisely, the quantity $\sup_{X\in\mathcal{X}}\vert f_v(X)\vert$ is bounded from above by a constant.  
This assumption  is easily satisfied as soon as the kernel $k_v$ is bounded on the compact set $\cX$, 
$$\sup_{X\in\mathcal{X}}\vert f_v(X)\vert \leq \sup_{X \in \cX} \sqrt{k_{v}(X_v, X_v)} \|f_v\|_{\cH_{v}}.$$ 
For a detailed discussion on this subject, we refer to the paper by \cite{Raskutti:2012:MRS:2503308.2188398}.
\item In the Gaussian regression model,
\begin{itemize}
\item \cite{Raskutti:2012:MRS:2503308.2188398} assumed that the unknown function $m$ has a sparse univariate decomposition, where each component in its decomposition lies in a RKHS. They obtained an estimator for $m$, based on a ridge group sparse type procedure. They established upper and lower bounds on the risk in the $L^2$-norm and upper bound on the risk in the empirical $L^2$-norm. 

\item \cite{huet:hal-01434895} assumed that the unknown function $m$ admits a Hoeffding decomposition involving the main effects and interactions. They obtained a RKHS ridge group sparse estimator of a meta-model that approximates the Hoeffding decomposition of $m$. They established upper bounds on the risk in the $L^2$-norm and the empirical $L^2$-norm. 
\end{itemize}
\cite{Raskutti:2012:MRS:2503308.2188398} and \cite{huet:hal-01434895} do not assume global boundedness condition. Instead, they assume that for all $v\in\mathcal{P}$ the functions $f_v$ are uniformly bounded. The proof of their results relies on the empirical Gaussian process methods such as Sudakov minoration \cite{Pisier1989} and concentration inequalities for Gaussian processes. 

As we are not in the Gaussian regression model, these methods could not be used in our work. We require new tools that we describe in details in the two next Sections.
\end{itemize}
\subsection{Sudakov minoration}\label{sec:sudakov} 
In the following Section, we recall the definition of the covering numbers, the statement of the classical Sudakov minoration, which is specific to the Gaussian process, and the generalized Sudakov minoration known also as the Sudakov minoration principal, which could be applied to some other processes. In Section \ref{ssudakov} we state the appropriate Sudakov type minoration to the process associated with the random variables that are distributed with density $\pi_\alpha\in\mathcal{D}$ (see Equation \eref{familypi}) in Corollary \ref{corosudakov}. 
\subsubsection{Introduction}\label{sintro}
Let $T$ be a set of square-integrable functions, i.e. $T\subset L^2$, and $\Vert.\Vert$ be the Euclidean norm. 
For any $\delta>0$, we denote by $C(\delta,T,\Vert .\Vert)$ the $\delta$-covering set of the metric space $(T,\Vert .\Vert)$:
\begin{align*}
C(\delta,T,\Vert .\Vert)=\Big\{f^1,...,f^N:\:\forall f\in T,\:\exists k\in\{1,...,N\}\text{ such that }\Vert f-f^k\Vert\leq \delta\Big\}.
\end{align*}
The $\delta$-covering number of $(T,\Vert .\Vert)$, denoted $N(\delta,T,\Vert .\Vert)$, is the cardinal of the smallest covering set. A proper covering restricts the covering to use only elements in the set $T$. It can be shown that the covering numbers and the proper covering numbers are related by the following inequality:
\begin{align}
\label{covering}
N(\delta,T,\Vert .\Vert)\leq N_{\text{proper}}(\delta,T,\Vert .\Vert)\leq N(\frac{\delta}{2},T,\Vert .\Vert).
\end{align} 
Consider a random variable $Z$ such that $E(Z^2) < \infty$, and consider an i.i.d. sequence $\{Z_i\}_{i=1}^n$ distributed like $Z$. To each $t =(t_1,...,t_n)$ of $T\subset L^2$ one can associate the process $V_t = \sum_{i=1}^n Z_i t_i$, $t\in T$.

In order to link the covering number on a class $T$, i.e. $N(\delta,T,\Vert .\Vert)$, to the expectation of the supremum of the process $V_t=\sum_{i=1}^n Z_i t_i$ in the Gaussian setting, the classical Sudakov minoration could be used (\cite{Pisier1989}):
\begin{align}
\label{csudv}
\frac{1}{K} \log N(\delta,T,\Vert .\Vert)\leq \Big(\frac{nE_{Z}\sup_{t\in T}\sum_{i=1}^nZ_i t_i}{\delta}\Big)^2.
\end{align}
When dealing with the processes $V_t=\sum_{i=1}^n Z_i t_i$, $t\in T$ associated with the random variables $\{Z_i\}_{i=1}^n$ that are not Gaussian, a generalized Sudakov minoration, known also as the Sudakov minoration principal, could be used to lower bound the value $E_Z\sup_{t\in T}\sum_{i=1}^n Z_i t_i$. Let us recall this inequality.

\begin{mydefinition}(Definition 1.1. in \cite{Latala2014SudakovtypeMF}) Let $Z=(Z_1,...,Z_n)$ be a random vector in $\mathbb{R}^n$. We say that $Z$ satisfies the $L_p$-Sudakov minoration principle with a constant $K'>0$, $SMP_p(K')$, if for any set $T\subset\mathbb{R}^n$ with $\vert T\vert>\exp(p)$ such that
\begin{align}
\label{smplata}
\Big(E_Z\sum_{t,s\in T}\vert \sum_{i=1}^n (t_i-s_i)Z_i\vert^p\Big)^{1/p}:=\Vert \sum_{i=1}^n (t_i-s_i)Z_i\Vert_p\geq \delta, \:\forall s,t\in T,\: s\neq t,
\end{align}
we have
\begin{align*}
K' \delta\leq E_Z \sup_{t,s\in T}\sum_{i=1}^n (s_i-t_i)Z_i.
\end{align*}
A random vector $Z$ satisfies the Sudakov minoration principle with a constant $K'$, $SMP(K')$, if it satisfies $SMP_p(K')$ for any $p\geq 1$.
\end{mydefinition}
If $\{Z_i\}_{i=1}^n$ are independent symmetric $\pm 1$ random variables or equivalently if the vector $Z=(Z_1,...,Z_n)$ is uniformly distributed on the cube $[-1,1]$ the Sudakov minoration principal with universal $K'$ was proven by \cite{talagrand1993}.

\cite{Latala2014SudakovtypeMF} proved the Sudakov minoration principal for the independent log-concave random variables. 
A measure on $\mathbb{R}^n$ with the full dimensional
support is log-concave if and only if it has a density of the form $\exp(-\phi(x))$, where $\phi:\mathbb{R}^n\rightarrow(-\infty,\infty]$ is convex (\cite{borell1974}). In the dependent setting the Sudakov minoration principal for the log-concave random variables was proven by \cite{Bednorz2014SomeRO}. 

As we are in the independent setting and the densities $\pi_\alpha\in\mathcal{D}$ (see Equation \eref{familypi}) are log-concave, the Sudakov minoration obtained by \cite{Latala2014SudakovtypeMF} holds in our context. 
However, we could not deduce from the result obtained by \cite{Latala2014SudakovtypeMF} the adapted Sudakov type minoration that leads to obtain the \textit{optimal} rate of convergence for our estimator. By \textit{optimal} we mean the same rate of convergence as in the Gaussian regression setting (see \cite{huet:hal-01434895}). 
This is the reason why we restricted ourselves to the densities $\pi_\alpha\in\mathcal{D}$ for which there exists a result given by \cite{10.2307/2374931}. 

In the next Section we provide in Corollary \ref{corosudakov} the appropriate Sudakov type minoration for the random variables that are distributed with density $\pi_\alpha\in\mathcal{D}$. This Corollary is a consequence of the result obtained by \cite{10.2307/2374931}.  
\subsubsection{Sudakov minoration for density $\pi_\alpha$}\label{ssudakov}
In this Section we state in Corollary \ref{corosudakov} the Sudakov minoration appropriate for the random variables that are distributed with density $\pi_\alpha\in\mathcal{D}$ (see Equation \eref{familypi}). This Corollary is a consequence of the Sudakov minoration stated in Theorem 3.1. in \cite{10.2307/2374931}. We start by introducing some notation that we need in the rest of this Section.

Let us denote by $\tilde{\alpha}$ the conjugate exponent of $\alpha$, i.e. $1/\alpha+1/\tilde{\alpha}=1.$ 
So, for all $\alpha>2$ we have $1<\tilde{\alpha}< 2$. 

We consider the sets $B_{\tilde{\alpha}}$ and $U_{\tilde{\alpha}}(u)$, $u\geq 0$ defined as follows:
\begin{align}
B_{\tilde{\alpha}}=\Big\{x\in\mathbb{R}^n:\:\sum_{k= 1}^n\vert x_k\vert^{\tilde{\alpha}}\leq 1\Big\},
\end{align}
and 
\begin{align}
\label{uset}
U_{\tilde{\alpha}}(u)=\Big\{x\in\mathbb{R}^n:\:\sum_{i=1}^n\eta_{\tilde{\alpha}}(x_i)\leq u,\: u\geq 0\Big\},
\end{align}
where 
$$\eta_{\tilde{\alpha}}(x_i)=x_i^2\mathbf{1}_{[-1,1]}(x_i)+\vert x_i\vert^{\tilde{\alpha}}\mathbf{1}_{(-\infty,-1]\cap[1,\infty)}(x_i).$$

For $T\subset L^2$ and $u\geq 0$, let $D(T,U_{\tilde{\alpha}}(u))$ be a covering set of translates of $T$ by $U_{\tilde{\alpha}}(u)$:
\begin{align*}
D(T,U_{\tilde{\alpha}}(u)) &=\Big\{f^1,...,f^N:\:\forall f\in T,\:\exists k\in\{1,...,N\}\text{ such that }f-f^k\in U_{\tilde{\alpha}}(u)\Big\},\\
&=\Big\{f^1,...,f^N:\:\forall f\in T,\:\exists k\in\{1,...,N\}\text{ such that }\sum_{i=1}^N\eta_{\tilde{\alpha}}(f_i-f_i^k)\leq u\Big\}.	
\end{align*}
We denote by $N(T,U_{\tilde{\alpha}}(u))$ the minimum number of translates of $U_{\tilde{\alpha}}(u)$ by elements of $T$ needed to cover $T$.

\begin{lemma}
\label{lem23tala}
For all $\tilde{\alpha}\leq 2$ and $u\geq0$, it is shown that (\cite{10.2307/2374931}): 
\begin{align}
\label{UBrelation}
U_{\tilde{\alpha}}(u)\subset (u^{1/2}B_2 + u^{1/\tilde{\alpha}}B_{\tilde{\alpha}} ).
\end{align}
\end{lemma}

\begin{rem}
\label{ubrelation}
If $\tilde{\alpha}\leq 2$ and $u\geq0$, then
$$U_{\tilde{\alpha}}(u)\subset 2\times \max (u^{1/2} ,u^{1/\tilde{\alpha}})B_2.$$
\end{rem}
The proof of Remark \ref{ubrelation} is given in Section \ref{proofubrelation} page \pageref{proofubrelation}.




\begin{theo} (Theorem 3.1. in \cite{10.2307/2374931})
\label{theotalag}
Let $Z=(Z_1,...,Z_n)$ be i.i.d. random variables distributed with density $\pi_\alpha\in\mathcal{D}$ defined in Equation \eref{familypi}, $U_{\tilde{\alpha}}(u)$, $u\geq 0$ be defined by \eref{uset} and $T\subset L^2$. Set 
\begin{align}
\label{msudakov}
M=E_{Z}\sup_{t\in T}\sum_{i=1}^nt_iZ_i,
\end{align}
then it is shown that:
\begin{align}
\label{talagSudakov}
N(T,U_{\tilde{\alpha}}(M))\leq \exp(KM),
\end{align}
where $K$ is a constant that depends on $\alpha$ only.
\end{theo}

\begin{rem}
\label{lemsudkvmax}
According to Theorem \ref{theotalag} and Remark \ref{ubrelation} for all $u\geq 0$ we have,
\begin{align}
\label{sudkvmax}
N(2\times \max(u^{1/2},u^{1/\tilde{\alpha}}),T,\Vert .\Vert)\leq N(T,U_{\tilde{\alpha}}(u))\leq\exp(Ku).
\end{align}
To be more precise, since $1<\tilde{\alpha}< 2$ we have
\begin{itemize}
\item[(i)] For $u\leq1$, $u^{1/\tilde{\alpha}}\leq u^{1/2}$ and $N(2u^{1/2},T,\Vert .\Vert)\leq\exp(Ku)$.
\item[(ii)] For $u\geq1$, $u^{1/\tilde{\alpha}}\geq u^{1/2}$ and $N(2u^{1/\tilde{\alpha}},T,\Vert .\Vert)\leq\exp(Ku)$.
\end{itemize}
\end{rem}

\begin{coro}
\label{corosudakov}
Under the same assumptions as for Theorem \ref{theotalag} we have for all $\delta>0$,
\begin{align*}
\frac{1}{K}\log N(\delta,T,\Vert .\Vert)\leq (\frac{2M}{\delta})^{\alpha}\mathbf{1}_{(0,2M]}(\delta)+(\frac{2M}{\delta})^2 \mathbf{1}_{[2M,\infty)}(\delta),
\end{align*}
which is exactly Equation \eref{sudtalagrand} with $M$ defined in Equation \eref{msudakov}.
\end{coro}
The proof of Corollary \ref{corosudakov} is given in Section \ref{proofcorosudakov} page \pageref{proofcorosudakov}.

\subsection{Concentration inequality}\label{sec:conineq}
We start this Section with a small introduction on the concentration inequalities context in Section \ref{introconc}, and we detail the concentration inequality used in our work in Section \ref{sconineq}.
\subsubsection{Introduction}\label{introconc}
Let $Z=(Z_1,...,Z_n)$ be a random vector in $\mathbb{R}^n$, and the function $\phi$ from $\mathbb{R}^n$ to $\mathbb{R}$ be convex and $1-$Lipschitz with respect to the Euclidean norm on $\mathbb{R}^n$, i.e.
$$\Vert \phi(Z)-\phi(Z')\Vert\leq\Vert Z-Z'\Vert,\: Z,Z'\in\mathbb{R}^n.$$

We are interested in the concentration inequalities of order two that provide bounds on how $\phi(Z)$ deviates from its expected value. More precisely, for  $P$ being the probability measure on $\mathbb{R}^n$, and for all $u\geq 0$,
\begin{equation}
\label{deviation}
P\Big(\vert \phi(Z)-E(\phi(Z)) \vert \geq u\ \Big)\leq C_1\exp\Big(-\frac{u^2}{C_2}\Big),
\end{equation}
where $C_1$, and $C_2$ are constants.

It was shown by \cite{LedouxTal:91} that, if $Z=(Z_1,...,Z_n)$ is a centered Gaussian random vector in $\mathbb{R}^n$, then:
\begin{align*}
P\Big(\vert \phi(Z)-E(\phi(Z))\vert\geq u\Big)\leq 4\exp\Big(-\frac{u^2}{2}\Big).
\end{align*}
This result could be proved using an inequality established by geometric arguments and an induction on the number of coordinates. 
 
After that, an alternative approach to some of Talagrand’s inequalities was proposed by \cite{PS1997160} based on the log-Sobolev inequalities. He showed that if the probability measure $P$ on $[0,1]^n$ satisfies the log-Sobolev inequality then it satisfies the concentration inequalities of the form (\ref{deviation}), i.e. the log-Sobolev inequality implies the deviation inequality. 

We say that the probability measure $P$ satisfies the log-Sobolev inequality for a class of functions $\Psi$ with loss function $R:\mathbb{R}^n\rightarrow[0,+\infty)$, if for every $\psi \in \Psi$ we have,
\begin{equation*}
\text{Ent}(\exp(\psi))\leq CE (R(\nabla\psi)\exp(\psi)),
\end{equation*}
where $\nabla\psi$ is the usual gradient of $\psi$, and $\text{Ent}(\exp(\psi))$ is the usual entropy of $\exp(\psi)$, i.e. 
$$\text{Ent}(\exp(\psi))=E (\psi \exp(\psi))-E (\exp(\psi))\log(E (\exp(\psi))).$$

This inequality was first introduced by \cite{10.2307/2373688} with $R(x) = \Vert x\Vert^2$, $x\in\mathbb{R}^n$ and $\Psi$ being the class of $\mathcal{C}^1$ functions. A lot of work has been done with different loss and class of functions, see for example \cite{Bobkov1997}, \cite{gentil:hal-00001609,gentil2007}. 

In the rest of this paper, we assume that $\Psi$ is the class of convex functions, and we consider only the quadratic loss $R(x)=\Vert x\Vert^2$, $x\in\mathbb{R}^n$. Therefore, the probability measure $P$ satisfies the convex log-Sobolev inequality if,  
\begin{equation}
\label{lsobolev}
E (\psi \exp(\psi))-E (\exp(\psi))\log(E (\exp(\psi)))\leq CE (\Vert\nabla\psi\Vert^2\exp(\psi)).
\end{equation}
\cite{adamarticle} found a sufficient condition for a class of probability distributions, denoted $\mathcal{M}(m,\rho^2)$ with $m>0$ and $\rho\geq0$, on the real line, to satisfy the convex log-Sobolev inequality. 
He deduced then the following concentration inequality which is satisfied for all probability distributions belonging to $\mathcal{M}(m,\rho^2)$:
\begin{align}
\label{concadam}
P\Big(\phi(Z)-E(\phi(Z))\geq u\Big)\leq \exp\Big(-\frac{u^2}{4C(m,\rho^2)}\Big).
\end{align}

We show in Lemma \ref{lipsc2} that the probability distributions associated with the densities $\pi_\alpha\in\mathcal{D}$ defined in Equation \eref{familypi} belong to $\mathcal{M}(m,\rho^2)$, and so they satisfy the convex log-Sobolev inequality.
As a consequence the concentration inequality \eref{concadam} holds for them. 

Recall that (see Section \ref{sketch} page \pageref{sketch}) we need concentration bounds for the
lower and upper tails of $\phi(Z)$, while the concentration inequality \eref{concadam} does not contain these two sides. 

\cite{articleshu} gave a sufficient and necessary condition for a probability measure on the real line to satisfy the convex log-Sobolev inequality. They obtained concentration bounds for the
lower and upper tails of convex functions of independent random variables
which satisfy the convex log-Sobolev inequality.

The result obtained by \cite{articleshu} allows us to state in Corollary \ref{lipschitz} the appropriate concentration inequality for the probability distributions associated with the densities $\pi_\alpha\in\mathcal{D}$.

\subsubsection{Concentration inequality for density $\pi_\alpha$}\label{sconineq}
In this Section we give the definition of the class of probability distributions $\mathcal{M}(m,\rho^2)$ and some of its properties. 
We show in Lemma \ref{lipsc2} that the probability distributions associated with the densities $\pi_\alpha\in\mathcal{D}$ (see Equation \eref{familypi}) belong to $\mathcal{M}(m,\rho^2)$, and so they satisfy the convex log-Sobolov inequality \eref{lsobolev}.
Finally, we state in Corollary \ref{lipschitz} the appropriate concentration inequality for our work which is a consequence of the concentration inequality stated in Corollary 1.7. of the paper by \cite{articleshu}.

\begin{mydefinition}
\label{mdef}
(Definition $4$ in \cite{adamarticle}) For $m>0$ and $\rho\geq 0$ let $\mathcal{M}(m,\rho^2)$ denote the class of probability distributions $\Pi$ on $\mathbb{R}$ for which
\begin{align*}
\upsilon^+(A)\leq \rho^2\Pi(A),
\end{align*}
for all sets $A$ of the form $A = [x,\infty)$, $x \geq m$ and
\begin{align*}
\upsilon^-(A)\leq \rho^2\Pi(A),
\end{align*}
for all sets $A$ of the form $A = (-\infty,-x]$, $x \geq m$, where $\upsilon^+$ is the measure on $[m,\infty)$
with density $ x\Pi([x,\infty))$ and $\upsilon^-$ is the measure on $(-\infty,-m]$ with density
$ -x\Pi((-\infty, x])$.
\end{mydefinition}

\begin{example}(Example page 5 in \cite{adamarticle})
\label{exadamczak}
The absolutely continuous distributions $\Pi$ that satisfy for $t\geq m$,
\begin{align}
\label{mfami}
\frac{d}{dt}\log \Pi ([t,\infty))\leq -\frac{t}{\rho^2} \quad \text{and}\quad 
\frac{d}{dt}\log \Pi ((-\infty,-t])\leq -\frac{t}{\rho^2}.
\end{align}
belong to $\mathcal{M}(m,\rho^2)$. In particular, if $\Pi$ has density of the form $\exp(-V (x))$ with $dV(x)/dx \geq x/\rho^2$ and $dV(-x)/dx \leq -x/\rho^2$ then $\Pi\in\mathcal{M} (1,\rho^2)$.
\end{example}

It is shown by \cite{adamarticle} that the probability distributions belonging to $\mathcal{M}(m,\rho^2)$ satisfy the convex log-Sobolov inequality \eref{lsobolev}. 
Let us denote by $\Pi_\alpha$ the probability distribution associated with the density $\pi_\alpha\in\mathcal{D}$ defined in Equation \eref{familypi}.
In the following Lemma we will show that $\bigotimes\Pi_\alpha$ satisfies the convex log-Sobolev inequality \eref{lsobolev}.

\begin{lemma}
\label{lipsc2} 
There exists some $m$ such that $\Pi_\alpha\in\mathcal{M}(m,\rho^2)$, and therefore $\bigotimes\Pi_\alpha$ satisfies the convex log-Sobolev inequality \eref{lsobolev}.
\end{lemma}
The proof of Lemma \ref{lipsc2} is given in Section \ref{prooflipsc2} page \pageref{prooflipsc2}.

As $\Pi_\alpha\in\mathcal{M}(m,\rho^2)$ and they satisfy the convex log-Sobolev inequality \eref{lsobolev}, so the concentration bound \eref{concadam} holds for them. Recall that (see Section \ref{sketch} page \pageref{sketch}), we need a concentration bound for the both upper and lower tails of a convex function of the random variables that are distributed as $\Pi_\alpha$. Therefore, the concentration bound \eref{concadam} is not sufficient for our work. We state in Corollary \ref{lipschitz} the appropriate concentration inequality for our work which is a consequence of the concentration inequality obtained by \cite{articleshu}. This result holds under a supplementary condition that we will state in the following Remark. 
\begin{rem} 
\label{expecfini}
Let $Z$ be a random variable distributed as $\Pi_\alpha$, then for every $s>0$ the quantity $E(\exp({s\vert Z\vert}))$ exists and is finite.
\end{rem}
The proof of Remark \ref{expecfini} is given in Section \ref{proofexpecfini} page \pageref{proofexpecfini}.

Note that, if $\alpha<2$ then $E(\exp({s\vert Z\vert}))\nless\infty$.
\begin{coro}
\label{lipschitz}
Let $Z=(Z_1, . . . ,Z_n)$ be i.i.d. random variables distributed as $\Pi_\alpha$.
Then there exists $A, B<\infty$ (depending only on $C$ in the log-Sobolev inequality \eref{lsobolev}), such that for any convex (or concave) function $\phi : \mathbb{R}^n\to\mathbb{R}$ which is $1-$Lipschitz (with respect to the Euclidean norm on $\mathbb{R}^n$) we have:
\begin{equation}
\label{conshu}
P\Big(\vert \phi(Z)-E(\phi(Z))\vert\geq u\Big)\leq 2B \exp\Big(-\frac{u^2}{8A}\Big),\: u\geq 0.
\end{equation}
\end{coro}
Corollary \ref{lipschitz} is a consequence of the concentration inequality shown by \cite{articleshu}:
\begin{equation}
\label{conshumed}
P\Big(\vert \phi(Z)-M(\phi(Z))\vert\geq u\Big)\leq B\exp\Big(-\frac{u^2}{A}\Big),\: u\geq 0,
\end{equation}
where $M$ is the median of $\phi(Z)$.

The proof of Corollary \ref{lipschitz} is given in Section \ref{prooflipschitz} page \pageref{prooflipschitz} and is based on the fact that the concentration inequalities around the mean and the median are equivalent up to a numerical constant (\cite{Milman:1986:ATF:21465}).

\section{Proof of Theorem \ref{oracle}}\label{prooforacle}
The proof  is based on four main  lemmas proved in Section~\ref{proofsintermediate}. In Section~\ref{intermediate} other lemmas used all along the proof are stated. 

Let us first establish inequalities that will be used in the following.
Let $f \in \cH$ and  $v\in S_f$ (see~\eref{sparsity}).

Using that for any $v \in S_f$, and any norm $\|\cdot
\|$ in $\cH_v$, $\| f_v\|-\| \widehat{f}_v\|
\leq \| f_v-\widehat{f}_v\|$ and that for any $v \notin
S_f$, $\| f_v\|=0$, 
we get,
\begin{align}
\sum_{v \in \mathcal{P}} \mu_v\| f_v\|_{\mathcal{H}_v}-\sum_{v \in \mathcal{P}} \mu_v\| \widehat{ f}_v\|_{\mathcal{H}_v}
\leq \sum_{v\in S_f} \mu_v\| f_v-\widehat{f}_v\|_{\mathcal{H}_v}-\sum_{v\notin S_f}
\mu_v\| \widehat{
  f}_v\|_{\mathcal{H}_v},\label{decomppen_H} 
\end{align}
and,
\begin{align}  
\sum_{v \in \mathcal{P}} \gamma_v\| f_v\|_{n}-\sum_{v \in \mathcal{P}}  \gamma_v\| \widehat{ f}_v\|_{n}
\leq \sum_{v\in S_f} \gamma_v\| f_v-\widehat{f}_v\|_{n}-\sum_{v\notin S_f}
\gamma_v\| \widehat{ f}_v\|_{n}. \label{decomppen_n}
\end{align}
Combining~\eref{decomppen_H},  and~\eref{decomppen_n},  to the fact that for any function  $f\in \mathcal{H}$, $\mathcal{L}(\widehat{f})\leq \mathcal{L}(f)$, we obtain,
\begin{equation*}
 \| m-\widehat{f}\|_n^2 \leq \| m-f\|_n^2+ B,
\end{equation*}
with
\begin{align}
B = 2V_{n,\varepsilon}\big(\widehat{f}-f\big)+
 \sum_{v\in S_f} [\mu_v\| \widehat{f}_v-f_v\|_{\mathcal{H}_v}+\gamma_v\| \widehat{f}_v-f_v\|_{n}]-\sum_{v\notin S_f}
[\mu_v\| \widehat{ f}_v\|_{\mathcal{H}_v}+  \gamma_v\| \widehat{f}_v\|_{n}].\label{B.eq}
\end{align}
If $\| m-f\|_n^2 \geq B$,  we immediately get the result
since in that case
\begin{equation*}
 \| m-\widehat{f}\|_n^2 \leq 
2 \| m-f\|_n^2 \leq 
2 \| m-f\|_n^2 + \sum_{v\in  S_f}\mu_v+\sum_{v\in S_f}\gamma^2_v.
\end{equation*}
If  $\| m-f\|_n^2 < B$, we get that
\begin{align}
\|
\widehat{f}-m\|_n^2 \leq  &2 B \label{base}\\
\leq  &4 \vert V_{n,\varepsilon}\big(\widehat{f}-f\big)\vert
+2\sum_{v\in S_f}[ \mu_v \|\widehat{f}_v-f_v \|_{\mathcal{H}_v}+\gamma_v \|\widehat{f}_v-f_v \|_{n}].
\label{base2.eq}
\end{align}
The control of the empirical process $\vert
V_{n,\varepsilon}\big(\widehat{f}-f\big)\vert$ is given by the
following lemma (proved in Section~\ref{ProofTau}, page~\pageref{ProofTau}).

 \begin{lemma} 
\label{Tau}
Let $V_{n,\varepsilon}$ be defined in \eref{procemp}.
For any $f$ in $\cF$, we
consider the event   $\mathcal{T}$ defined as
\begin{align}
\label{evtTau}\mathcal{T}=\left\lbrace \forall f \in \cF, \forall
  v \in \mathcal{P} , \vert
  V_{n,\varepsilon}\big(\widehat{f}_v-f_v\big) 
\vert \leq  
\kappa \lambda_{n,v}^2\| \widehat{f}_v-f_v\|_{\mathcal{H}_v}
+ 
\kappa \lambda_{n,v}\| \widehat{f}_v-f_v\|_{n}\right\rbrace,
\end{align}
where $\lambda_{n,v}$ is defined in Equation~\eref{lambda.eq} and where $\kappa= 10 + 4\Delta$. Then, for some positive constants $c_1, c_2$,
\begin{align}
\label{conrtolproba}
\PXe \left(\mathcal{T}\right)\geq 1-c_1\sum_{v \in \mathcal{P}}\exp(-nc_2 \lambda_{n,v}^2).
\end{align}
\end{lemma}
Conditioning on  $\mathcal{T}$, Inequality~\eref{base2.eq} becomes
\begin{align*}
 \| \widehat{f}-m\|_n^2 \leq & 
4 \kappa \sum_{v \in \mathcal{P}} [ \lambda_{n,v}^2 \|
\widehat{f}_v-f_v\|_{\mathcal{H}_v} 
+ \lambda_{n,v}\| \widehat{f}_v-f_v\|_{n}]+\\ 
&2\sum_{v\in S_f} [\mu_v\|
  \widehat{f}_v-f_v\|_{\mathcal{H}_v}+ \gamma_v\| \widehat{f}_v-f_v\|_{n}],
\end{align*}
which may be decomposed as follows
\begin{align*}
\| \widehat{f}-m\|_n^2 \leq&   
\sum_{v\in S_f} [
4\kappa \lambda_{n,v}^2+2\mu_v] \|
\widehat{f}_v-f_v\|_{\mathcal{H}_v} + \sum_{v\in
  S_f}[
4\kappa \lambda_{n,v}+2\gamma_v]\|
\widehat{f}_v-f_v\|_{n} +
\\& 4\sum_{v\notin S_f} 
\kappa \lambda_{n,v}^2 \| \widehat{f}_v-f_v\|_{\mathcal{H}_v} +
4\sum_{v\notin S_f} 
\kappa \lambda_{n,v}\| \widehat{f}_v-f_v\|_{n}.
\end{align*}
If we choose $C_{1}\geq \kappa$ in Theorem~\ref{oracle}, then $\kappa \lambda_{n,v}^{2} \leq \mu_{v}$
and $\kappa \lambda_{n,v}\leq  \gamma_{v}$ and the previous inequality becomes 
\begin{align}
  \|  \widehat{f}-m\|_n^2  \leq &6 \sum_{v\in S_f}[\mu_v\| \widehat{f}_v-f_v\|_{\mathcal{H}_v}+\gamma_v\| \widehat{f}_v-f_v\|_{n}]+\nonumber\\
&4\sum_{v\notin S_f}[\mu_v \| \widehat{f}_v\|_{\mathcal{H}_v}+\gamma_v
\| \widehat{f}_v\|_{n}]. \label{base3}
\end{align}
Next we use the decomposability property of the penalty expressed in the
following lemma (proved in Section~\ref{Prooflemme2} page~\pageref{Prooflemme2}).

\begin{lemma} 
\label{lemme2}
For any $f \in \mathcal{F}$, under the assumptions of Theorem \ref{oracle},
conditionally on $\mathcal{T}$ (see~\eref{evtTau}), we have: 
\begin{equation}
\label{lemme2eq}
 \sum_{v\notin S_f} \mu_v \|\widehat{
   f}_v\|_{\mathcal{H}_v} +\sum_{v\notin S_{f}}
 \gamma_v\| \widehat{ f}_v\|_{n} \leq  
3 \sum_{v\in S_f} \mu_v\|
\widehat{f}_v-f_v\|_{\mathcal{H}_v}+
3 \sum_{v\in S_f} \gamma_v\| \widehat{f}_v-f_v\|_{n}.
\end{equation}
\end{lemma}
Hence, by combining \eref{base3}  and Lemma~\ref{lemme2} we obtain 
\begin{align*}
 \|  \widehat{f}-m\|_n^2\leq  18 
\sum_{v\in S_f}\big[\mu_v \| \widehat{f}_v-f_v\|_{\mathcal{H}_v} +\gamma_v \| \widehat{f}_v-f_v
\|_{n
}\big].\end{align*}
For each $v$, $\|
\widehat{f}_v-f_v\|_{\mathcal{H}_v} \leq 2$ (because the functions  $\widehat{f}_v$ et $f_v$ belong to the
class $\cF$, see~\eref{calF.eq}), and 
consequently, for some constant $C$,
\begin{equation}
 \label{borne}
 \| \widehat{f}-m\|_n^2  \leq C\Big\{
 \sum_{v\in  S_f}\mu_v+
 \sum_{v\in S_f}\gamma_v \| \widehat{f}_v-f_v
\|_{n}\Big\}.
\end{equation}
To finish the proof it remains to compare the two quantities $\sum_{v\in S_f} \| \widehat{f}_v-f_v
\|_{n}^2$ and  $\| \sum_{v\in S_f} \widehat{f}_v-f_v
\|_{n}^2$. For that purpose we show that 
$\| \sum_{v\in S_f}\widehat{f}_v-f_v\|_n$ is less than $ \|
\sum_{v\in S_f}\widehat{f}_v-f_v\|_2^2$ plus an
additive term coming from concentration results (see the Lemma given
below). Next, thanks to the orthogonality of the spaces $\cH_{v}$
with respect to $L^{2}(P_X,\mathcal{X})$, 
$ \|\sum_{v\in  S_f}\widehat{f}_v-f_v\|_2^2=\sum_{v \in  S_f} 
\|\widehat{f}_{v}-f_{v}\|_2^2$. To
conclude,  it remains to
consider several cases, according to the rankings of 
$\|\sum_{v\in  S_f}\widehat{f}_v-f_v\|_2^2$ and
$\|\sum_{v\in  S_f}\widehat{f}_v-f_v\|_n^2$.
This is the subject of the following lemma whose proof is given in Section~\ref{Proofnorm2normn}, page~\pageref{Proofnorm2normn}.

\begin{lemma} 
\label{norm2normn}
For $f \in \mathcal{H}$, let $\mathcal{A}$ be the  event 
\begin{equation}
\label{A}
\mathcal{A}=\Big\{ \forall f \in \mathcal{F}, \forall v \in \mathcal{P}, \;
 \|\widehat{f}_v-{f}_v\|_n \leq  2\| \widehat{f}_v-{f}_v\|_2
+\gamma_v\Big\}.
\end{equation}
Then, for some positive constant $c_{2}$,
\begin{equation*}
P_{X,\varepsilon} ( \mathcal{A} ) \geq 1-\sum_{v\in\mathcal{P}} \exp(-n c_{2}\gamma_v^2).
\end{equation*}
\end{lemma}

On the set $\cA$, Inequality \eref{borne}  provides that, for all $K>0$
\begin{align}
\frac{1}{C} \| \widehat{f}-m\|_n^2 &\leq 
\sum_{v\in  S_f} [\mu_v + 2 \gamma_v \| \widehat{f}_v-{f}_v\|_2 + \gamma_v^2
], \nonumber\\
&\leq \sum_{v\in  S_f} [\mu_v + (1+K)\gamma_v^{2}  +\frac{1}{K}\| \widehat{f}_v-{f}_v\|^{2}_2
], \label{intermediaire1}\\
&\leq \sum_{v\in  S_f} [\mu_v + (1+K)\gamma_v^{2} ] +\frac{1}{K}\sum_{v \in \mathcal{P}}\| \widehat{f}_v-{f}_v\|^{2}_2, \nonumber\\
\label{intermediaire2}
 &\leq  \sum_{v\in  S_f} [\mu_v + (1+K)\gamma_v^{2}]  
  + \frac{1}{K}\| \sum_{v \in \mathcal{P}}\widehat{f}_v-{f}_v \|^{2}_2.
\end{align}
Inequality~\eref{intermediaire1} uses the inequality  $2ab \leq
\frac{1}{K} a^{2}+ K b^{2}$ for all positive $K$, and
Inequality~\eref{intermediaire2} uses the orthogonality with respect
to $L^{2}(\PX)$.

In the following we have to consider several cases, according to the
rankings of
$ \|\sum_{v\in \mathcal{P}}\widehat{f}_v-{f}_v\|_2$ and $\| \sum_{v\in \mathcal{P}}\widehat{f}_v-{f}_v\|_n$. 
More precisely, we consider two following cases:
\begin{itemize}
\item[\underline{Case 1:}] If $ \|  \sum_{v\in \mathcal{P}}\widehat{f}_v-{f}_v\|_2\leq \| \sum_{v\in \mathcal{P}}\widehat{f}_v-{f}_v\|_n$.
\item[\underline{Case 2:}] If $\| \sum_{v\in \mathcal{P}}\widehat{f}_v-{f}_v\|_2\geq \| \sum_{v\in \mathcal{P}}\widehat{f}_v-{f}_v\|_n $.
\end{itemize}
\medskip
\noindent
\label{3Cases}
\underline{Case 1:}
From \eref{intermediaire2}, for any $f \in  \mathcal{H}$, we get
\begin{align*}
\frac{1}{C}\| \widehat{f}-m\|_n^2 
 \leq  
 \sum_{v\in  S_f} [\mu_v+ (1+K)
\gamma_v^2 ]+
  \frac{1}{K}\| \widehat{f}-{f}\|_n^2.
\end{align*}
Hence, using that for all $K'>0$, 
\begin{equation}
\label{trick2}\|\widehat{f}-f\|_{n}^{2}\leq (1+K') \|\widehat{f} -m\|_{n}^2+
(1+\frac{1}{K'})\|f-m\|_{n}^2,
\end{equation} 
we obtain for a suitable choice of $K'$, say $1+K' < K/C$, that, for some positive constant $C'$,
\begin{align*}
\| \widehat{f}-m\|_n^2 
  \leq  C'\Big\{
\| f-m\|_n^2+ 
\sum_{v\in  S_f}\mu_v+ 
 \sum_{v\in S_f}\gamma_v ^2\Big\}.
\end{align*}
This shows the result in Case 1.\\ 
\medskip
\noindent
\underline{Case 2:} This case is solved by applying the following Lemma (proved in
Section~\ref{prooflemcompnormes3},
page~\pageref{prooflemcompnormes3}), which states that with high
probability, $\|\widehat{f} -f \|_{2} \leq \sqrt{2} \|\widehat{f} -f \|_{n}$.

\begin{lemma} 
\label{lemcompnormes3}
Let $f=\sum_v f_v \in \cF$ with support $S_{f}$, $\lambda_{n,v}$ be
defined by~\eref{lambda.eq},  and let
$\mathcal{G}(f)$ be the class of functions written as $g=\sum_{v\in
  \mathcal{P}} g_v$, such that $\| g_v\|_{\mathcal{H}_v}\leq 2$
satisfying for all $f \in \cF$
\begin{eqnarray*}
\mbox{\bf C1}  && \sum_{v\in \mathcal{P}} \mu_v\| g_v \|_{\mathcal{H}_v} +\sum_{v\in \mathcal{P}} \gamma_v\| g_v\|_{n}\leq  
4\sum_{v\in S_f} \mu_v\| g_v\|_{\mathcal{H}_v}+
4\sum_{v\in S_f} \gamma_v\| g_v\|_{n} \\
\mbox{\bf C2}  && \sum_{v\in S_f}\gamma_v\| g_v\|_n\leq 2\sum_{v\in
  S_f}\gamma_v\| g_v\|_2+\sum_{v\in S_f} \gamma_v^2\\
\mbox{\bf C3}  && \| g\|_n\leq \| g \|_2\\
\end{eqnarray*}
Then the event 
\begin{align*}
 \Big\{  \| g\|_n^2 \geq \frac{\| g\|^2_2}{2}\Big\},
\end{align*}
have probability greater than $1-c_1\exp(-n c_3 \sum_{v\in S_f}\lambda_{n,v}^2)$ for some constants $c_1$ and $c_3$.
\end{lemma} 

If $f$ is such that $|S_{f}|=0$, then Condition {\bf C1} is not satisfied except if $g_{v}=0$ for all $v  \in \mathcal{P}$. Because we will apply Lemma~\ref{lemcompnormes3} to $g_{v}= \widehat{f}_{v} - f_{v}$, this event has probability $0$.
If $f$ is such that $|S_{f}| \geq 1$, then Condition {\bf C1} is satisfied:

from Equation \eref{lemme2eq} in Lemma \ref{lemme2} we have, 
\begin{align*}
&\sum_{v\notin S_f} \mu_v \|\widehat{f}_v\|_{\mathcal{H}_v}+\sum_{v\in S_f} \mu_v\|
\widehat{f}_v-f_v\|_{\mathcal{H}_v} +\sum_{v\notin S_{f}} \gamma_v\| \widehat{ f}_v\|_{n}+ \sum_{v\in S_f} \gamma_v\| \widehat{f}_v-f_v\|_{n}\\
&\leq  
3 \sum_{v\in S_f} \mu_v\|\widehat{f}_v-f_v\|_{\mathcal{H}_v}+\sum_{v\in S_f} \mu_v\|
\widehat{f}_v-f_v\|_{\mathcal{H}_v}+
3 \sum_{v\in S_f} \gamma_v\| \widehat{f}_v-f_v\|_{n}+\sum_{v\in S_f} \gamma_v\| \widehat{f}_v-f_v\|_{n},\\
&\Leftrightarrow \sum_{v\in \mathcal{P}} \mu_v\|\widehat{f}_v-f_v\|_{\mathcal{H}_v}+\sum_{v\in \mathcal{P}} \gamma_v\| \widehat{f}_v-f_v\|_{n}\leq 4\sum_{v\in S_f} \mu_v\|
\widehat{f}_v-f_v\|_{\mathcal{H}_v}+4 \sum_{v\in S_f} \gamma_v\| \widehat{f}_v-f_v\|_{n}.
\end{align*} 
Moreover, Assumption $n \lambda_{n,v}^{2} \geq -C_{2}\log(\lambda_{n,v})$ implies that  
$$\lambda_{n,v}= K_{n,v}/\sqrt{n} \mbox{ with } K_{n,v}\rightarrow \infty.$$
Then,
$$\exp(-n c_{3}\sum_{v \in S_{f}} \lambda_{n,v}^2) \leq \exp(-c_{3}\vert S_{f}\vert\min_{v \in \mathcal{P}} K_{n,v}^2),$$ 
and the event 
\begin{equation}
\label{C}
\mathcal{C}=\Big\{ \forall f \in \mathcal{F}, \mbox{ such that } 
 g=\sum_{v\in \mathcal{P}}(\widehat{f}_v-f_v )\in \mathcal{G}(f),
\mbox{ and }
 \| g\|_n^2 \geq \frac{\| g\|^2_2}{2} \Big\}
\end{equation}
 has probability greater than $1-\eta/3$ for some $0<\eta<1$.

Conditioning on the events $\mathcal{T}$ and  $\mathcal{A}$
(defined by \eref{evtTau} and \eref{A}), $\sum_{v\in \mathcal{P}}
(\widehat{f}_v-{f}_v)$ belongs to the set $\mathcal{G}(f)$.
According to \eref{intermediaire2}, we conclude in the same way as in the first
case.

Finally, it remains to quantify $P_{X,\varepsilon}(\mathcal{T} \cap \mathcal{A}  \cap \mathcal{C})$. Following
Lemma~\ref{Tau}, and Lemma~\ref{norm2normn}, $\mathcal{T}$, respectively $\mathcal{A}$,
has probability greater than $1 - c_1\sum_{v
  \in \mathcal{P}}\exp(-nc_2 \lambda_{n,v}^2)$, respectively
$1-\sum_{v
  \in \mathcal{P}} \exp(-n \gamma_v^2)$. Each of these probabilities is greater than
$1-\eta/3$ thanks to the assumption $n \lambda_{n,v}^{2} \geq -C_{2}
\log \lambda_{n,v}$.

\hfill $\Box$

\subsection{Intermediate Lemmas}
\label{intermediate}
\begin{lemma} 
\label{lemcomplex}
If $\EXe$ denotes the expectation with respect to the distribution of $(X, \varepsilon)$,  we have for all $t > 0 $,
\begin{equation*}
\EXe  W_{n,2,v}(t)
 \leq  Q_{n,v}(t).
\end{equation*}
\end{lemma}
Its proof is given in Section~\ref{Prooflemcomplex} page~\pageref{Prooflemcomplex}.

\begin{lemma} 
\label{lemcompnormes1}
Let  $b > 0$ and let $\cG(t)$ be the following class of functions:
\begin{equation}
\cG(t)=\Big\{ g_v \in \cH_{v}, \|g_v\|_{\cH_v}\leq 2, \|g_v\|_{2}\leq t,
  \|g_v\|_{\infty}\leq b \Big\}.
\label{calG.eq}
\end{equation}
Let $\Omega_{v,t}$ be the event defined as
\begin{align}
\label{Omega}
\Omega_{v,t}=\Big \{
\sup_{g_v \in \cG(t)} \{
\vert \| g_v\|_2- \| g_v\|_n \vert \} \leq \frac{bt}{2}\Big\}.
\end{align}
Then for any $t\geq \nu_{n,v}$, the event $\Omega_{v,t}$ has
probability greater than $ 1-\exp(-c_2 n t^2)$, 
for some positive constant $c_{2}$.
\end{lemma}
Its proof is given in Section~\ref{Prooflemcompnormes1},
page~\pageref{Prooflemcompnormes1}.

\begin{lemma} 
\label{lemcompnormes2}
For any function $g_v \in \mathcal{H}_v$
satisfying $\| g_v\|_{\mathcal{H}_v}\leq 2$,   $\| g_v\|_\infty \leq b$ and 
 $\| g_v\|_2\geq t$, for all $t\geq
   \nu_{n,v}$ and $b \geq 1$, the event
\begin{equation*}
(1 - \frac{b}{2}) \| g_v\|_2 \leq \| g_v\|_n \leq 
(1 + \frac{b}{2}) \| g_v\|_2
\end{equation*}
has probabilty greater than $1-\exp(-c_2 n t^2)$ for some positive
constant $c_{2}$.
\end{lemma}
Its proof is given in Section~\ref{Prooflemcompnormes2},
page~\pageref{Prooflemcompnormes2}.

\begin{lemma} 
\label{concentration1}
If $\Ee$ denotes the expectation with
respect to the distribution of $\varepsilon$,  we have
\begin{align}
\label{concentnn}\PXe \Big( 
\vert W_{n,n,v}(t) - E_\varepsilon \big( W_{n,n,v}(t)\big) \vert\geq
\delta t   \Big) \leq 2B\exp(- \frac{n\delta^2}{8A}).
\end{align}
\end{lemma}
Its proof is given in Section~\ref{Proofconcentration1},
page~\pageref{Proofconcentration1}.

\begin{lemma}
\label{concentration2}
Conditionally on the space $\Omega_{v,t}$ defined by \eref{Omega}, we
have the following inequalities:
\begin{align}
\label{concentn2}
\PXe \Big( \vert W_{n,2,v}(t) - E_\varepsilon\big( W_{n,2,v}(t)\big) 
\vert \geq \delta t  \Big)  \leq  2B\exp(-
  \frac{n\delta^2}{32A}), 
\end{align}
\begin{align}
\label{concentration}
\PX \Big(  \Ee W_{n,2,v}( t) - \EXe
\big( W_{n,2,v}(t)\big)  \geq  x \Big) \leq 
\exp( -\frac{n x^2}{ Q_{n,v}(t)} ).
\end{align}
\end{lemma}
Its proof is given in Section~\ref{Proofconcentration2},
page~\pageref{Proofconcentration2}.

 \begin{lemma}
  \label{Case1} Let $\lambda_{n,v}$ be defined at
  Equation~\eref{lambda.eq}, $\Delta$ at Equation~\eref{nu} and
  $\kappa= 10+4\Delta$. Conditionally
on the space  $\Omega_{v,\lambda_{n,v}}$ defined at Equation~\eref{Omega}, for some positive constants $c_1, c_{2}$, 
 with probability greater than $1-c_{1}\exp(-c_{2} n\lambda_{n,v}^2)$,
we have 
\begin{equation}
 \label{etap} W_{n,n,v}( \lambda_{n,v}) \leq \kappa \lambda_{n,v}^{2}\;
 \mbox{ and } \; \Ee W_{n,n,v}( \lambda_{n,v}) \leq \kappa \lambda_{n,v}^{2}.
\end{equation}
\end{lemma}
Its proof is given in Section~\ref{ProofCase1}, page~\pageref{ProofCase1}.

\subsection{Proof of lemma \ref{Tau} to \ref{lemcompnormes3}}\label{proofsintermediate}
\subsubsection{Proof of lemma \ref{Tau}}
\label{ProofTau}
For $f \in \cF$ and $v \in \cP$, let $g_{v} =
\widehat{f}_v-f_v$. Note that $\|g_{v}\|_{\cH_{v}} \leq 2$. Let us show that 
\begin{equation}
\label{but}
\vert V_{n,\varepsilon}(g_v)\vert 
\leq 
\kappa\Big(\lambda_{n,v}^2\| g_v\|_{\mathcal{H}_v}+\lambda_{n,v}\| g_v\|_{n}\Big).
\end{equation}
We start by writing that
\begin{equation}
| V_{n,\varepsilon}(g_v) |=
\| g_v\|_{\mathcal{H}_v}\Big| V_{n,\varepsilon}\Big(\frac{g_v}{\| g_v\|_{\mathcal{H}_v}}\Big) \Big|
\leq \| g_v\|_{\mathcal{H}_v} W_{n,n,v}\Big(\frac{\| g_v\| _{n}}{ \|
  g_v\|_{\mathcal{H}_v}}\Big).
\label{baseW}
\end{equation}
Consider the two following cases:
\begin{itemize}
\item[\underline{Case A:}] $\| g_v\|_n \leq \lambda_{n,v} \|
  g_v\|_{\mathcal{H}_v}$,
\item[\underline{Case B:}] $\| g_v\|_n > \lambda_{n,v} \| g_v\|_{\mathcal{H}_v}$.
\end{itemize}

\medskip
\noindent
\underline{Case A}: \label{CaseA}
Since $\| g_v\| _{n}\leq \lambda_{n,v}\| g_v\|_{\mathcal{H}_v} $, we have
$$W_{n,n,v}\Big(\frac{\| g_v\| _{n}}{ \|
  g_v\|_{\mathcal{H}_v}}\Big)\leq W_{n,n,v}( \lambda_{n,v}).$$ We
then apply Lemma \ref{Case1}, page~\pageref{Case1}, and conclude that \eref{but} holds in Case A for each  $v\in \mathcal{P}$ since,
 with high probability
\begin{align}
\label{cas1}
\vert V_{n,\varepsilon}(g_v)\vert 
\leq 
\kappa\lambda_{n,v}^2\| g_v\|_{\mathcal{H}_v}\leq 
\kappa\lambda_{n,v}^2\| g_v\|_{\mathcal{H}_v}+
\kappa\lambda_{n,v}\| g_v\|_{n}.
  \end{align}

\medskip
\noindent\underline{Case B}: Consider now the case $\| g_v\| _{n}>\lambda_{n,v}
\| g_v\| _{\mathcal{H}_v}$ and let us show that for any $v\in \mathcal{P}$,
\begin{align*}
W_{n,n,v}(\frac{\| g_v\|_n}{\| g_v\|_{\mathcal{H}_v}})
 \leq 
\kappa\lambda_{n,v}\| g_v\|_{n}
 .
\end{align*}
Let $r_v$ be a deterministic number such that $r_v>\lambda_{n,v}$. Our first step relies on the study of the process
 $W_{n,n,v}(r_v),$ for  $r_v>\lambda_{n,v}$.
In that case we state two results:
\begin{itemize}
\item[{\bf R1}]  For any deterministic $r_v\geq \lambda_{n,v}$, with probability greater than $1-c_1 \exp(-c_2 n \lambda_{n,v}^2)$,
\begin{equation}\label{int43}
W_{n,n,v}(r_v) \leq \kappa r_v \lambda_{n,v}.
\end{equation}

\item[{\bf R2}] Inequality~\eref{int43} continues to hold
  for random $r_v$ of the form $$r_v=\frac{\| g_v\|_n}{\|
    g_v\|_{\mathcal{H}_v}}.$$
\end{itemize}
Combining these two points implies that, with probability greater than $1-c_1 \exp(-c_2 n \lambda_{n,v}^2)$,
\begin{equation*}
 \| g_v\|_{\mathcal{H}_v}W_{n,n,v}\Big(\frac{\| g_v\|_n}{\|
  g_v\|_{\mathcal{H}_v}}\Big) \leq
\kappa \| g_v\|_n \lambda_{n,v}.
\end{equation*}

Consequently, in Case B, according to \eref{baseW}, for each $v$,
Inequality~\eref{but} holds because
\begin{equation*}
\vert V_{n,\varepsilon}(g_v)\vert 
\leq 
\kappa \| g_v\|_n \lambda_{n,v}
\leq 
\kappa \lambda_{n,v}^2\| g_v\|_{\mathcal{H}_v}+ 
\kappa \lambda_{n,v}\| g_v\|_{n}.
\end{equation*}
 This ends up the proof of Lemma \ref{Tau}.

\medskip
\paragraph{\underline{Proof of {\bf R1}}}

From Lemma \ref{concentration1}, page~\pageref{concentration1} with $t=r_v$ and $\delta=\lambda_{n,v}$, we get that with  probability
 greater than $1-2B\exp(-n\lambda_{n,v}^2/8A)$,
 \begin{equation}
 W_{n,n,v}(r_v)\leq E_\varepsilon(W_{n,n,v}(r_v))+r_v\lambda_{n,v}
 \end{equation}
Next we prove that for some positive  $r_v$, with probability greater than $1-\exp(-nc\lambda_{n,v}^2)$, we have
 \begin{align}
\label{intcas2}\Ee (W_{n,n,v}( r_v))\leq
\kappa  r_v \lambda_{n,v}.
\end{align}
Let $\widehat{\nu}_{n,v}$ be defined 
as the smallest solution of $\Ee ( W_{n,n,v}(t))\leq \kappa t^2$.
For $W_{n,n,v}$, defined by \eref{defWntilde}, we write
\begin{align*}
\Ee (W_{n,n,v}(r_v))&= \frac{r_v}{\widehat{\nu}_{n,v}}\Ee  \sup \Big\{\vert V_{n,\varepsilon}(g_v)\vert , \;
\| g_v\|_{\mathcal{H}_v}\leq 2(\frac{\widehat{\nu}_{n,v}}{r_v}) , \;
\| g_v\|_{n}\leq \widehat{\nu}_{n,v}\Big\}.
\end{align*}
Besides, Lemma~\ref{Case1} stated that on the event $\Omega_{v,  \lambda_{n,v}}$, $\Ee (W_{n,n,v}( \lambda_{n,v}) )\leq \kappa\lambda_{n,v}^{2}$. It follows from the definition of
$\widehat{\nu}_{n,v}$, and Lemma~\ref{lemcompnormes1}, that $\widehat{\nu}_{n,v} \leq \lambda_{n,v}$ for all $v \in\mathcal{P}$ with probability greater than $1 - \exp(-n c_{2} \sum_{v  \in \cP} \lambda_{n,v}^{2})$.  Consequently, for any deterministic $r_v$ such that $r_v\geq \lambda_{n,v}$, we have
\begin{align*}
\widehat{\nu}_{n,v}\leq \lambda_{n,v}\leq r_v\Leftrightarrow \frac{\widehat{\nu}_{n,v}}{r_v}\leq 1,
\end{align*}
and so,
\begin{align*}
\Ee (W_{n,n,v}(r_v))&= \frac{r_v}{\widehat{\nu}_{n,v}}\Ee  \sup \Big\{\vert V_{n,\varepsilon}(g_v)\vert , \;
\| g_v\|_{\mathcal{H}_v}\leq 2 , \;
\| g_v\|_{n}\leq \widehat{\nu}_{n,v}\Big\},\\
&\leq \frac{r_v}{\widehat{\nu}_{n,v}}  \Ee ( W_{n,n,v}(\widehat{\nu}_{n,v}))\leq \frac{r_v}{\widehat{\nu}_{n,v}}\kappa \widehat{\nu}_{n,v} ^2=\kappa  r_v \widehat{\nu}_{n,v}\leq \kappa  r_v \lambda_{n,v}.
\end{align*}
\medskip
\paragraph{\underline{Proof of {\bf R2}}}

Let us  prove  {\bf R2} by using a peeling-type argument.  Our aim is to prove that \eref{int43} holds for any $r_v$ of the form
$$r_v=\frac{\| g_v\|_n}{\| g_v\|_{\mathcal{H}_v}}.$$ 
Since $\|g_v\|_\infty/\| g_v\|_{\mathcal{H}_v} \leq 1$, we have
$\|g_v\|_n/\| g_v\|_{\mathcal{H}_v} \leq 1$. 
We thus restrict ourselves  to $r_v$ satisfying $r_v= \| g_v\|_n/\|
g_v\|_{\mathcal{H}_v}$ with $\| g_v\|_n/\| g_v\|_{\mathcal{H}_v} \in (\lambda_{n,v},1]$.

We start by splitting the interval $(\lambda_{n,v},1]$ into $M$ disjoint intervals such that 
$$(\lambda_{n,v},1]=\cup_{k=1}^M(2^{k-1}\lambda_{n,v},2^k\lambda_{n,v}],$$ 
for some $M$ that will be chosen later. Consider the event $\cD^{c}$ defined as follows:
\begin{align*}
\mathcal{D}^c=\Big\{ \exists v \in \mathcal{P} \mbox{ and } \exists  \overline{g}_v, \mbox{ such that }   
\vert V_{n,\varepsilon}(\overline{g}_v)\vert 
\geq  
\kappa \lambda_{n,v}\| \overline{g}_v \|_{n}, 
\mbox{ with } \frac{\| \overline{g}_v\|_n}{\| \overline{g}_v\|_{\mathcal{H}_v}} \in (\lambda_{n,v},1]\Big\}.
\end{align*}

We prove that, for some positive constants $c_1, c_2$,
$$P (\mathcal{D}^c) \leq c_1\exp(-c_2 n\lambda_{n,v}^2).$$

For $\overline{g}_v \in \cD^{c}$, let $\overline{k}$ be the integer in $\{1,\cdots,M\}$, such that
$$2^{\overline{k}-1}\lambda_{n,v} \leq   \frac{\| \overline{g}_v\|_n}{\| \overline{g}_v\|_{\mathcal{H}_{v}}}  \leq 2^{\overline{k}}\lambda_{n,v}.$$
This $\overline{k}$ satisfies
\begin{align*}
\| \overline{g}_v\|_{\mathcal{H}_v}W_{n,n,v}\Big(  2^{\overline{k}}\lambda_{n,v}\Big)\geq
\| \overline{g}_v\|_{\mathcal{H}_v}W_{n,n,v}\Big(  \frac{\| g_v\|_n}{\| g_v\|_{\mathcal{H}_v}}\Big)
\geq \vert V_{n,\varepsilon}(\overline{g}_v)\vert
\geq
\kappa \lambda_{n,v}\| \overline{g}_v\|_{n}.
\end{align*}
Therefore, we get
$$W_{n,n,v}(2^{\overline{k}}\lambda_{n,v}) \geq
\kappa \lambda_{n,v}\frac{\|
  \overline{g}_v\|_n}{\| \overline{g}_v\|_{\mathcal{H}_v}}\geq
\kappa \lambda_{n,v}^2 2^{\overline{k}-1}
\geq 
\kappa \frac{\lambda_{n,v}}{2} 2^{\overline{k}}\lambda_{n,v}.$$
By taking $r_v=2^{\overline{k}}\lambda_{n,v}$ in \eref{int43}, we have
$$\mathcal{P} \Big( W_{n,n,v}(2^{\overline{k}}\lambda_{n,v}) 
\geq 
\kappa \frac{\lambda_{n,v}}{2}
2^{\overline{k}}\lambda_{n,v}\Big)\leq c_1 \exp(-c_2n \lambda_{n,v}^2).$$

Now let us write $\mathcal{D}^c$ as
follows:
\begin{align*}
 \mathcal{D}^c =\bigcup_{k=1}^M\Big\{\; 
\exists v\mbox{ and } \exists \; \overline{g}_v \mbox{ such that }
\vert V_{n,\varepsilon}(\overline{g}_v)\vert 
\geq 
\kappa \lambda_{n,v}\| \overline{g}_v\|_n,
\mbox{ with } \frac{\| \overline{g}_v\|_n}
{\| \overline{g}_v \|_{\mathcal{H}}} 
\in (2^{k-1}\lambda_{n,v},2^k\lambda_{n,v}]\Big\}.
\end{align*}
The set
$\mathcal{D}^c$ has probability smaller than  $c_1 M \exp(-c_2n \lambda_{n,v}^2)$.
If we choose $M $ such that $\log M \leq (c_{2}/2) n
\lambda_{n,v}^{2}$, then the probability of the set $\cT$ is  greater than 
\begin{equation*}
1-\sum_{v \in\cP}  c_{1} \exp(- \frac{c_{2}}{2} n \lambda_{n,v}^{2}).
\end{equation*}

It follows that \textbf{R2} is proved which ends up the proof of Lemma \ref{Tau}.

\hfill $\Box$ 

\subsubsection{Proof of lemma \ref{lemme2}}
\label{Prooflemme2}
Starting from \eref{base} with $B$ defined by Equation~\eref{B.eq}, we write
\begin{align*}
 \frac12\| \widehat{f}-m\|_n^2 \leq &
2\vert V_{n,\varepsilon}(\widehat{f}-f)\vert +  \sum_{v\in S_f}[\mu_v\|\widehat{f}_v-f_v\|_{\mathcal{H}_v}+\gamma_v \|\widehat{f}_v-f_v\|_{n}]-\\
& \sum_{v\notin S_f}[\mu_v\| \widehat{ f}_v\|_{\mathcal{H}_v} +\gamma_v \| \widehat{f}_v\|_{n}].
\end{align*}
On the event $\mathcal{T}$ defined in \eref{evtTau} we have
\begin{align*}
\frac12 \| \widehat{f}-m\|_n^2 \leq  &2 \kappa \sum_{v\in \mathcal{P}} \lambda_{n,v}^2 \|
\widehat{f}_v-f_v\|_{\mathcal{H}_v}+2\kappa\sum_{v\in \mathcal{P}} \lambda_{n,v}\| \widehat{f}_v-f_v\|_{n}+\\
&\sum_{v\in S_f}[\mu_v  \| \widehat{f}_v-f_v        \|_{\mathcal{H}_v}+\gamma_v \| \widehat{f}_v-f_v       \|_{n}]
- \sum_{v\notin S_f}
[\mu_v\| \widehat{ f}_v\|_{\mathcal{H}_v}+\gamma_v\| \widehat{ f}_v\|_{n}].
\end{align*}
Rearranging the terms we obtain that
\begin{align*}
\frac12\| \widehat{f}-m\|_n^2 \leq 
&\sum_{v\in S_f}(2 
\kappa\lambda_{n,v}^2+\mu_v)\|
\widehat{f}_v-f_v\|_{\mathcal{H}_v}+\sum_{v\in S_f}(2
\kappa\lambda_{n,v}+\gamma_v)\| \widehat{f}_v-f_v\|_{n}+\\
&\sum_{v\notin S_{f}}
(2 
\kappa\lambda_{n,v}^2 -\mu_v)\| \widehat{ f}_v\|_{\mathcal{H}_v}+\sum_{v\notin S_{f}}
(2 
\kappa\lambda_{n,v} -\gamma_v)\| \widehat{ f}_v\|_{n}.
\end{align*}
Now, thanks to Assumption~\eref{condmu} with $C_{1} \geq \kappa$
  we have  $\kappa\lambda_{n,v}^2
  \leq \mu_v$  and $2 \kappa\lambda_{n,v}
  \leq\gamma_v$ and Lemma \ref{lemme2} is shown  since
\begin{align*}
0\leq \frac12\| \widehat{f}-m\|_n^2 \leq &  
3\sum_{v\in S_f}
\mu_v\| \widehat{f}_v-f_v\|_{\mathcal{H}_v}
+3\sum_{v\in S_f}
\| \widehat{f}_v-f_v\|_{n}-
\\ &\sum_{v\notin S_{f}} \mu_v\| \widehat{ f}_v\|_{\mathcal{H}_v} -\sum_{v\notin S_{f}} \gamma_v\| \widehat{ f}_v\|_{n}.
\end{align*}

\hfill $\Box$

\subsubsection{Proof of lemma \ref{norm2normn}}
\label{Proofnorm2normn}
Let us consider the following two cases:
\begin{itemize}
\item $\|\widehat{f}_v-{f}_v\|_2 \leq \gamma_v$. We apply
  Lemma~\ref{lemcompnormes1} (page~\pageref{lemcompnormes1}) to the
  function  $g_{v}=\widehat{f}_v-{f}_v$. It satisfies $g_{v}\in
  \cG(\gamma_v)$ with $b=2$ (recall that $\|\cdot\|_{\infty} \leq
  \|\cdot\|_{\cH_{v}}$). 
  Moreover, $\gamma_{v} \geq C_{1}
  \lambda_{n,v}\geq C_{1} \nu_{nv} \geq \nu_{n,v}$ as soon as $C_1 \geq
  1$.

It follows that, for some positive $c_{2}$, with probability greater than $1-\exp(-n c_{2}\gamma_{v}^{2})$,
\begin{equation*}
\| \widehat{f}_v-{f}_v\|_n\leq
  \| \widehat{f}_v-{f}_v\|_2+\gamma_v.
\end{equation*}

\item $\|\widehat{f}_v-{f}_v\|_2 \geq \gamma_v$. We apply
  Lemma~\ref{lemcompnormes2} (page~\pageref{lemcompnormes2}) to the
  function  $g_{v}=\widehat{f}_v-{f}_v$ with $b=2$. 
It follows that, for
some positive $c_{2}$, with probability greater than $1-\exp(-n c_{2}\gamma_{v}^{2})$,
\begin{equation*}
\| \widehat{f}_v-{f}_v\|_n\leq
  2 \| \widehat{f}_v-{f}_v\|_2.
\end{equation*}
\end{itemize}

\hfill $\Box$

\subsubsection{Proof of lemma \ref{lemcompnormes3}}
\label{prooflemcompnormes3}
Throughout the proof, we make use of the quantity $d_n$ defined as follows:  

For $\beta<1/\alpha$ and some constant $\eta'$, 
\begin{equation}
\label{dn.eq}
d_n^2\geq\eta'n^{\alpha\beta-1}.
\end{equation}
Let $\mathcal{G}(f)$ and $\mathcal{G}'(f)$ be the following sets:
\begin{align*}
 \mathcal{G}(f) =&
\Big\{ g=\sum_{v\in \mathcal{P}}g_v, \mbox{ satisfying }
 \|g_{v}\|_{\cH_v} \leq 2, \mbox{ and Conditions \textbf{C1, C2, C3}
 }\Big\},\\
\mathcal{G}'(f)  = & \Big\{ g\in \mathcal{G}(f),
\mbox{ such that } \| g\|_2=d_n\Big\}.
\end{align*}
In order to prove this lemma we consider two cases: if $ \| \sum_{v\in \mathcal{P}}\widehat{f}_v-{f}_v\|_2\geq d_n$, and if $ \| \sum_{v\in \mathcal{P}}\widehat{f}_v-{f}_v\|_2\leq d_n$.

First, we suppose that $\| \sum_{v\in \mathcal{P}}\widehat{f}_v-{f}_v\|_2\geq d_n$, and we consider the two events $\cB$ and $\cB'$ defined as follows:
\begin{equation*}
\mathcal{B}=\Big\{\forall h \in \mathcal{G}, \; \| h\|_n^2\geq \frac{\| h\|_2^2}{2} , \mbox{ and } \| h\|_2\geq d_n\Big\},
\end{equation*} 
and
\begin{equation}
\label{eventBprim}
\mathcal{B}'=\Big\{\forall h \in \mathcal{G}', \; \| h\|^{2}_n \geq
\frac{d_n^{2}}{2}  \Big\}.
\end{equation}
If $h \in \cB'$ , then $h \in \cG$, $\|h\|_{2} =d_n$ and $\| h\|_n^2\geq d_n^2/2$. It follows that  $\| h\|_n^2\geq\| h\|_2^2/2$ and $\|h\|_{2} \geq d_n$. We just showed that the event $\mathcal{B}'$ is included into the event $\mathcal{B}$. So, this case is proved if the event $\mathcal{B'}$ holds with high probability. Consider 
$$Z_n( \mathcal{G}')  =\sup_{g\in \mathcal{G}'}\Big\{d_n^2-\| g\|_n^2\Big\}.$$ 
We show that the event $Z_n( \mathcal{G}')\leq d_n^2/2$ has probability greater than $1-c_1  \exp(-nc_3 d_n^2)$.

Consider a $d_n/8$-covering  of $(\mathcal{G}', \|\cdot \|_{n})$. So that, for all $g$ in $\mathcal{G}'$ there exists $g^k$ such that $$\| g-g^k\|_n \leq \frac{d_n}{8}.$$  The associated proper covering number  is:
\begin{equation}
 N_{\rm{pr}}=N_{\rm{pr}}(\frac{d_n}{8},\mathcal{G}', \|\cdot \|_{n} ).
\label{CovN.eq}
\end{equation}
Now, for all $g \in \mathcal{G}'$, we write:
\begin{equation}
\label{t1t2}
 d_n^2-\| g\|_n^2=T_1+T_2,
\end{equation}
with $T_1=\| g^k\|_n^2-\| g\|_n^2$ and $T_2=d_n^2-\| g^k\|_n^2$.
The proof  is splitted into four steps:
\begin{enumerate}
\item[\underline{Step 1}] The first step consists in showing that 
\begin{align}
\label{etape1}
T_1=\| g^k\|_n^2-\|  g\|_n^2\leq \frac{d_n^2}{4}.
\end{align}

\item[\underline{Step 2}] The second step consists in proving that, for $N_{\rm{pr}}$ given at Equation~\eref{CovN.eq} and for some constant $C$,
\begin{align*}
\PX \Big( \max_{k\in\{1,\cdots,N_{\rm{pr}}\}}[d_n^2-\| g^k\|_n^2]\geq
  \frac{d_n^2}{4}   \Big)
\leq \exp\Big(\log{N_{\rm{pr}}}- C n d_n^2 \Big).
\end{align*}

\item[\underline{Step 3}]
The third  step concerns the control of  $N_{\rm{pr}}$. Let $\sigma^2_\alpha$ be the variance of a random variable distributed with density $\pi_\alpha\in\mathcal{D}$ (see Equation \eref{familypi}), then for some $K >0$,  
\begin{align*}
\frac{1}{K} \log N_{\rm{pr}} \leq&\Big(32\sigma_\alpha\sqrt{n}(\Ee\sup _{g\in\cG'}\vert V_{n,\varepsilon}(g)\vert)/d_n\Big)^{\alpha}\mathbf{1}_{(0,32\sigma_\alpha\sqrt{n} E_{\varepsilon}\sup _{g\in\cG'}\vert V_{n,\varepsilon}(g)\vert]}(d_n)+\\
&\mathbf{1}_{[32\sigma_\alpha\sqrt{n} E_{\varepsilon}\sup _{g\in\cG'}\vert V_{n,\varepsilon}(g)\vert,\infty)}(d_n).
\end{align*}

\item[\underline{Step 4}]
The last step consists  in bounding from above the Gaussian
complexity. For some $ \kappa>0$
\begin{equation*}
E_\varepsilon\sup_{g\in \mathcal{G}'} \sum_{v \in \cP}\vert V_{n,\varepsilon}(g_v)\vert 
 \leq \frac{4\kappa}{C_{1}} \Big\{ \sum_{v\in S_f} (2\mu_v   +
          \gamma_v^2) + 
2(\sum_{v\in S_f}\gamma_{v}^2)^{\frac{1}{2}}
 d_n\Big\},
\end{equation*}
\end{enumerate}
Let us conclude the proof of the lemma before proving these four steps.

Putting together Steps $3$ and $4$ we have:

If $d_n\in[32\sigma_\alpha\sqrt{n} E_{\varepsilon}\sup _{g\in\cG'}\vert V_{n,\varepsilon}(g)\vert,\infty)$, then 
\begin{equation*}
\frac{1}{K}  
\log N_{\rm{pr}}\Big(\frac{d_n}{8}, \mathcal{G}^\prime, \parallel
.\parallel_n\Big)
\leq 1.
\end{equation*}
Thanks to Step $2$, 
\begin{equation*}
\PX \Big( T_2\geq\frac{d_n^2}{4}\Big)\leq\PX \Big( \max_{k\in\{1,\cdots,N_{\rm{pr}}\}}[d_n^2-\| g^k\|_n^2]\geq
  \frac{d_n^2}{4}   \Big)
\leq K \exp\Big(- C n d_n^2 \Big),
\end{equation*}
and, therefore
\begin{equation}
\label{casless}
\PX \Big( Z_n(\mathcal{G}')\leq \frac{d_n^2}{2} \Big)
\leq K \exp\Big(- C n d_n^2 \Big).
\end{equation}
If $d_n\in(0,32\sigma_\alpha\sqrt{n} E_{\varepsilon}\sup _{g\in\cG'}\vert V_{n,\varepsilon}(g)\vert]$, then
\begin{align*}
\frac{1}{K} \log N_{\rm{pr}}\Big(\frac{d_n}{8}, \mathcal{G}^\prime, \parallel
.\parallel_n\Big)
 &\leq(32\sigma_\alpha)^\alpha n^{\frac{\alpha}{2}}\Big(\frac{\Ee\sup _{g\in\cG'}\vert V_{n,\varepsilon}(g)\vert}{d_n}\Big)^{\alpha},\\
&\leq(32\sigma_\alpha)^\alpha n^{\frac{\alpha}{2}}\Big(\frac{4\kappa}{C_{1} d_n} 
 ( \sum_{v\in S_f} (2\mu_v   +  \gamma_v^2) + 
2( \sum_{v\in S_f}\gamma_{v}^2)^{\frac{1}{2}}d_n)\Big)^{\alpha},\\
&\leq \Big(\frac{128 \kappa\sigma_\alpha}{C_{1}}\Big)^{\alpha}n^{\frac{\alpha}{2}}\Big(
\frac{\sum_{v\in S_f} (2\mu_v   +  \gamma_v^2) }{d_n}
+ 2 ( \sum_{v\in S_f}\gamma_{v}^2)^{\frac{1}{2}} \Big)^{\alpha}.
\end{align*}

We have to show that $\log N_{\rm{pr}} - C n d_n^{2} \leq -c_{3}
n d_n^{2}$ or equivalently that $\log N_{\rm{pr}} \leq \widetilde{C}n
d_n^{2}$, where $\widetilde{C} = C-c_{3}$.

Let $A =K (128 \kappa\sigma_\alpha/C_{1})^{\alpha}$. We have,
\begin{eqnarray*}
 \log N_{\rm{pr}}  \leq \widetilde{C}n d_n^{2}
& \Leftrightarrow & 
A n^{\frac{\alpha}{2}} \Big(
\frac{\sum_{v\in S_f} (2\mu_v   +  \gamma_v^2) }{d_n}
+ 2  (\sum_{v\in S_f}\gamma_{v}^2)^{\frac{1}{2}} \Big)^{\alpha}
\leq \widetilde{C}n d_n^{2}, \\
&  \Leftrightarrow &  
\Big(
\frac{\sum_{v\in S_f} (2\mu_v   +  \gamma_v^2) }{d_n}
+ 2  (\sum_{v\in S_f}\gamma_{v}^2)^{\frac{1}{2}} \Big)^{\alpha}
\leq \frac{\widetilde{C}}{A} n^{1-\frac{\alpha}{2}} d_n^{2}, 
 \\
&  \Leftrightarrow &  
\frac{\sum_{v\in S_f} (2\mu_v   +  \gamma_v^2) }{d_n}
+ 2  (\sum_{v\in S_f}\gamma_{v}^2)^{\frac{1}{2}} 
\leq (\frac{\widetilde{C}}{A})^{\frac{1}{\alpha}} n^{\frac{1}{\alpha}-\frac{1}{2}} d_n^{\frac{2}{\alpha}}. 
\end{eqnarray*}
Because $\gamma_{v} = C_{1}\lambda_{n,v}$ and $\mu_{v} = C_{1}
\lambda_{n,v}^{2}$, 
\begin{eqnarray*}
 \log N_{\rm{pr}}  \leq \widetilde{C}n d_n^{2}
& \Leftrightarrow & 
C_{1}(2  + C_{1}) \frac{\sum_{v\in S_f} \lambda_{n,v}^{2} }{d_n}
+ 2 C_{1} (\sum_{v\in S_f}\lambda_{n,v}^2)^{\frac{1}{2}} 
\leq (\frac{\widetilde{C}}{A})^{\frac{1}{\alpha}} n^{\frac{1}{\alpha}-\frac{1}{2}} d_n^{\frac{2}{\alpha}}. 
\end{eqnarray*}
Considering the first term in the left hand side, let 
\begin{equation*}
B =  \frac{1}{2}\times \frac{1}{C_{1}(2+C_{1})}(\frac{\widetilde{C}}{A})^{\frac{1}{\alpha}}, 
\end{equation*}
then
\begin{equation*}
\frac{\sum_{v\in S_f} \lambda_{n,v}^{2} }{d_n} \leq  B n^{\frac{1}{\alpha}-\frac{1}{2}}
d_n^{\frac{2}{\alpha}} \Leftrightarrow 
d_n^{2} \geq B^{-\frac{2\alpha}{2+\alpha}} \Big(\sum_{v\in S_f}
  \lambda_{n,v}^{2}\Big)^{\frac{2\alpha}{2+\alpha}} n^{\frac{\alpha-2}{\alpha+2}}.
\end{equation*}
As $\sum_{v\in S_f}\lambda_{n,v}^2 \leq C_3
n^{2\beta-1}$ (see Equation \eref{hypdf.eq}), we get 
\begin{equation*}
 B^{-\frac{2\alpha}{2+\alpha}} \Big(\sum_{v\in S_f}
  \lambda_{n,v}^{2}\Big)^{\frac{2\alpha}{2+\alpha}}
n^{\frac{\alpha-2}{\alpha+2}}
\leq  (\frac{B}{C_3})^{-\frac{2\alpha}{2+\alpha}} 
n^{\frac{4\alpha\beta}{2+\alpha}-1}.
\end{equation*}
Therefore, the inequality 
\begin{equation*}
C_{1}(2  + C_{1}) \frac{\sum_{v\in S_f} \lambda_{n,v}^{2} }{d_n}
\leq \frac{1}{2}(\frac{\widetilde{C}}{A})^{\frac{1}{\alpha}} n^{\frac{1}{\alpha}-\frac{1}{2}} d_n^{\frac{2}{\alpha}}, 
\end{equation*}
will be satisfied if 
\begin{equation*}
d_n^{2} \geq (\frac{C_3}{B})^{\frac{2\alpha}{\alpha+2}} n^{\frac{4\alpha\beta}{\alpha+2}-1} .
\end{equation*}

For the second term, let 
\begin{equation*}
B' = \frac{1}{2}\times\frac{1}{2C_{1}}  (\frac{\widetilde{C}}{A})^{\frac{1}{\alpha}},
\end{equation*}
then
\begin{equation*}
 \big(\sum_{v\in S_f}\lambda_{n,v}^2\Big)^{\frac{1}{2}} \leq B' n^{\frac{1}{\alpha}-\frac{1}{2}}
d_n^{\frac{2}{\alpha}}\Leftrightarrow 
d_n^{2} \geq B^{'- \alpha}\Big(\sum_{v\in S_f}
  \lambda_{n,v}^{2}\Big)^{\frac{\alpha}{2}} n^{\frac{\alpha-2}{2}}.
\end{equation*}
As $\sum_{v\in S_f}\lambda_{n,v}^2 \leq C_3
n^{2\beta-1}$ (see Equation \eref{hypdf.eq}), then 
\begin{equation*}
B^{'- \alpha}\Big(\sum_{v\in S_f}
  \lambda_{n,v}^{2}\Big)^{\frac{\alpha}{2}} n^{\frac{\alpha-2}{2}}
\leq \Big(\frac{C_3}{B^{'2}}\Big)^{\frac{\alpha}{2}}
n^{\alpha\beta-1}.
\end{equation*}
Therefore the inequality 
\begin{equation*}
2 C_{1} \Big(\sum_{v\in S_f}\lambda_{n,v}^2\Big)^{\frac{1}{2}} 
\leq \frac{1}{2}\Big(\frac{\widetilde{C}}{A}\Big)^{\frac{1}{\alpha}} n^{\frac{1}{\alpha}-\frac{1}{2}} d_n^{\frac{2}{\alpha}},  
\end{equation*}
will be satisfied if 
$$d_n^{2} \geq \Big(\frac{C_3}{B^{'2}}\Big)^{\frac{\alpha}{2}}n^{\alpha\beta-1}.$$

As $\alpha > 2$, $4\alpha\beta/(\alpha+2) < \alpha\beta$. Therefore,
there exists a constant $\eta'$, take for example 
$$\eta'=\max\Big(
(\frac{C_3}{B^{'2}})^{\frac{\alpha}{2}},
(\frac{C_3}{B})^{\frac{2\alpha}{\alpha+2}}\Big),$$ 
such that if $d_n^{2} \geq \eta' n^{\alpha\beta-1}$, then $\log N_{{\rm pr}} \leq \widetilde{C} n d_n^{2}$, and Step $2$ states that
\begin{equation*}
 \PX\left( T_{2} \geq \frac{d_n^{2}}{4}\right) \leq \PX \left( \max_{k\in\{1,\cdots,N_{\rm{pr}}\}}[d_n^2-\| g^k\|_n^2]\geq
  \frac{d_n^2}{4}   \right)
\leq \exp\Big(- c_{3} n d_n^{2}\Big).
\end{equation*}
Now, we have
\begin{equation}
\label{casgr}
 \PX \left(Z_n( \mathcal{G}')\leq \frac{d_n^2}{2}\right) =
\PX \left(\max_{g^1,\cdots,g^N} [  d_n^2-\| g^k\|_n^2  ]\geq \frac{d_n^2}{4}\right)\leq \exp\Big(-c_{3} nd_n^2\Big).
\end{equation} 
Finally, we obtain for $c_1=\max(K,1)$ and $c_3\leq C$ (see Equations \eref{casless} and \eref{casgr}):
\begin{equation*}
 \PX \left(Z_n( \mathcal{G}')\leq \frac{d_n^2}{2}\right) \leq c_1 \exp\Big(-c_{3} nd_n^2\Big).
\end{equation*} 
Moreover, for $n$ large enough, we have $\sum_{v \in S_{f}}
  \lambda_{n,v}^{2}\leq d_n^{2} \leq 1$ (see Equations \eref{hypdf.eq} and \eref{dn.eq}), and
  $$1- c_1\exp\Big(-c_{3}n
d_n^{2}\Big) \geq 1- c_1\exp\Big(- c_{3}n \sum_{v \in S_{f}} \lambda_{n,v}^{2}\Big).$$
Therefore,
\begin{equation*}
 \PX \left(Z_n( \mathcal{G}')\leq \frac{d_n^2}{2}\right) \leq c_1 \exp\Big(-c_{3} n\sum_{v \in S_{f}} \lambda_{n,v}^{2}\Big).
\end{equation*} 

Before proving the Steps $1$ to $4$ let us solve the second case: if $\| \sum_{v\in \mathcal{P}}\widehat{f}_v-{f}_v\|_2\leq d_n$ then we consider the event $\mathcal{B}''$ defined as follows:
\begin{equation*}
\mathcal{B}''=\Big\{\forall h \in \mathcal{G}, \; \| h\|^{2}_n \geq
\frac{\|h\|_2^{2}}{2} , \mbox{ and } \|h\|_{2} \leq d_n \Big\}.
\end{equation*}
We have that the event $\cB'$ defined in Equation \eref{eventBprim} is included in $\cB''$ and the same proof as in the first case applies.

\medskip

\paragraph{\underline{Proofs of Steps $1$ to $4$}}

\medskip

The proofs of Step $1$ and Step $2$ are strictly the same as in the Gaussian case.
More precisely

\medskip
\noindent
\underline{Proof of Step 1: }\label{ProofS1}

It is easy to see that,
\begin{align*}
T_1=\Vert g^k\Vert_n^2-\Vert g\Vert_n^2&=\frac{1}{n}\sum_{i=1}^n[(g^k(X_i))^2-(g(X_i))^2]\\
&=\frac{1}{n}\sum_{i=1}^n[g^k(X_i)-g(X_i)][g^k(X_i)+g(X_i)]\\
&\leq \Vert g^k-g\Vert_n\Big(\frac{1}{n}\sum_{i=1}^n[g^k(X_i)+g(X_i)]^2\Big)^{\frac{1}{2}}
\end{align*}
where in the inequality above we used Cauchy Schwarz inequality.
Using the inequality $(a+b)^2\leq 2a^2+2b^2$, $g \in \cG'$, and the property that $g$ satisfies
Condition {\bf C3}, we get
\begin{align*}
\frac{1}{n}\sum_{i=1}^n [g^k(X_i)+ g(X_i)]^2 \leq 2\| g^k\|_n^2+2\| g\|_n^2
\leq 4 d_n^2.
\end{align*}
Besides, the covering set is constructed such that $\| g^k-g\|_n \leq d_n/8$. It follows that Step 1 is proved.

\medskip
\noindent
\underline{Proof of Step 2:}

We prove that for some constant $C$, 
$$\PX \Big( T_2 \geq \frac{d_n^2}{4}\Big)\leq
\PX \Big(\max_{1 \leq k \leq N_{\rm{pr}}} [  d_n^2-\|
  g^k\|_n^2  ]\geq \frac{d_n^2}{4}\Big)
 \leq
\exp\left(\log{N_{\rm{pr}}}- C n d_n^{2} \right).$$
As $g^{k} \in \cG'$, $d_n = \| g^k\|_2$. Then 
$$\max_{1\leq k\leq N_{\rm{pr}}} [d_n^2-\| g^k\|_n^2 ]=\max_{1\leq k\leq N_{\rm{pr}}}[\| g^k\|_2^2-\| g^k\|_n^2].$$ 
Applying Theorem 3.5. in \cite{chung2006} with $X=\sum_i(g^k(X_i))^2$, for all positive $\lambda$ we have:
 $$\PX \Big(\sum_{i=1}^n [g^k(X_{i}))^2\leq n\mathbb{E} (g^k(X_{i})]^2   -\lambda\Big)\leq \exp\Big(   -\frac{\lambda^2}{2n \mathbb{E}(g^k(X))^4}\Big),$$
 or equivalently,
 $$\PX \Big(\| g^k\|_2^2-\| g^k\|_n^2\geq \frac{\lambda}{n}\Big)\leq \exp\Big(   -\frac{\lambda^2}{2n \mathbb{E}(g^k(X))^4}\Big).$$
 Taking $\lambda=nd_n^2/4$ and using that $\| g^{k}\|_2^2=d_n^2$ we get
 $$\PX \Big(d_n^2-\| g^k\|_n^2\geq \frac{d_n^2}{4}\Big)\leq \exp\Big(-\frac{nd_n^4}{32\mathbb{E}(g^k(X))^4}\Big).$$

 It follows that 
 \begin{align}\label{step21}
 \PX \Big(\max_{1\leq k\leq N_{\rm{pr}} } [d_n^2-\| g^k\|_n^2] \geq \frac{d_n^4}{4}\Big) &\leq \sum_{k=1}^{N_{\rm{pr}}}
\exp\Big(-\frac{nd_n^4}{32\mathbb{E}(g^k(X))^4}\Big) \nonumber\\
&\leq \exp\Big(\log{N_{\rm{pr}}}-\frac{n d_n^4}{32\max_{k}\mathbb{E}(g^k(X))^4}\Big).
 \end{align}

Moreover, $g \in\cH$, so $g=\sum_{v  \in  \cP} g_{v}$, where the functions $g_{v}$ are centered
and orthogonal in $L^{2}(\PX)$. Therefore $\E (g(X))^4$ is the sum of the following terms:
\begin{align*}
A_1&=\sum_{v\in \mathcal{P}} \EX g^4_v(X_v), \\
A_2 &=\begin{pmatrix}4\\2\end{pmatrix} \sum_{v \not= v'}
\EX g^2_v(X_v)g^2_{v'}(X_{v'}),\\
A_3&= 
\begin{pmatrix}4\\3\end{pmatrix} \sum_{v_1 \not= v_2\not= v_3}
\EX g^2_{v_1}(X_{v_1})
  g_{v_2}(X_{v_2}) g_{v_3}(X_{v_3}),\\
 A_4&=\begin{pmatrix} 4\\3\end{pmatrix}
\sum_{v_1 \not= v_2}
\EX  g^3_{v_1}(X_{v_1})g_{v_2}(X_{v_2}),\\
A_{5}&=\begin{pmatrix}4\\1\end{pmatrix} \sum_{v_1 \not= v_2\not= v_3\not=v_4}
\EX g_{v_1}(X_{v_1})g_{v_2}(X_{v_2})g_{v_3}(X_{v_3})g_{v_4}(X_{v_4}).
\end{align*}
Using the Cauchy Schwartz inequality and the fact that $\|
g_v\|_\infty\leq \| g_v\|_{\mathcal{H}_v}\leq 2$, and
$\|g\|_{2}=d_n$ (because $g\in \cG'$), we get that $A_{1}$
is proportional to $d_n^{2}$, $A_{2}, A_{3}, A_{5}$ to $d_n^{4}$,
and $A_{4}$ to $d_n^{3}$. For example, 
\begin{equation*}
A_1=\sum_{v\in \mathcal{P}} \EX g^{4}_v(X_v)\leq \|
g\|_\infty^2 \sum_{v\in \mathcal{P}}
\|g_{v}\|_{2}^{2}= \|
g\|_\infty^2 
\|\sum_{v\in \mathcal{P}}g_{v}\|_{2}^{2} \leq 4 d_n^2.
\end{equation*}
After calculation of the terms $A_i$, since $d_n^{2}$ is assumed to be smaller than one,
we get that:
\begin{align}
\label{step22}
\max_k \EX (g^k(X))^{4}\leq Cd_n^2(1+O(d_n^2)).
\end{align}
Step $2$ is proved  by combining \eref{step21} and \eref{step22}.

We now focus on Step $3$ and Step $4$:

\medskip
\noindent
\underline{Proof of Step $3$:} 

Let $N_{\rm{pr}}$ be defined at Equation~\eref{CovN.eq}. We  prove that
\begin{align*}
\frac{1}{K}\log N_{\rm{pr}}\Big(\frac{d_n}{8}, \mathcal{G}^\prime, \parallel .\parallel_n\Big)
\leq &\Big(32\sigma_\alpha\sqrt{n}(\Ee\sup _{g\in\cG'}\vert V_{n,\varepsilon}(g)\vert)/d_n\Big)^{\alpha}\mathbf{1}_{(0,32\sigma_\alpha\sqrt{n} E_{\varepsilon}\sup _{g\in\cG'}\vert V_{n,\varepsilon}(g)\vert]}(d_n)+\\
&\mathbf{1}_{[32\sigma_\alpha\sqrt{n} E_{\varepsilon}\sup _{g\in\cG'}\vert V_{n,\varepsilon}(g)\vert,\infty)}(d_n).
\end{align*}

We start from Equation \eref{covering} and write that:
$$\log N_{\rm{pr}}\Big(\frac{d_n}{8},\mathcal{G}',  
\|\cdot \|_{n} \Big)\leq \log N\Big(\frac{d_n}{16},\mathcal{G}',  
\|\cdot \|_{n} \Big).$$ 
Next, we use Corollary \ref{corosudakov}: 

Let $Z=(Z_1,...,Z_n)$ be i.i.d. random variables distributed with density $\pi_\alpha\in\mathcal{D}$ defined in Equation \eref{familypi} with $\sqrt{\mbox{var}(Z_i)}=\sigma_\alpha$.
 Set $T=\mathcal{G}'$, $\delta =\sqrt{n}d_n/16$ and $M=n\times E_{Z}\sup _{g\in\cG'}\vert V_{n,Z}(g)\vert$, then for all $\alpha\geq 2$ we have, 
\begin{align*}
\log N\Big(\frac{d_n}{16},\mathcal{G}',\Vert .\Vert_n\Big)=&\log N\Big(\frac{\sqrt{n}d_n}{16},\mathcal{G}',\Vert .\Vert\Big),\\
\leq &K\Big(\frac{32nE_{Z}\sup _{g\in\cG'}\vert V_{n,Z}(g)\vert}{\sqrt{n}d_n}\Big)^{\alpha}\mathbf{1}_{(0,2n\times E_{Z}\sup _{g\in\cG'}\vert V_{n,Z}(g)\vert]}(\frac{\sqrt{n}d_n}{16})+\\
&K\Big(\frac{32nE_{Z}\sup _{g\in\cG'}\vert V_{n,Z}(g)\vert}{\sqrt{n}d_n}\Big)^{2}\mathbf{1}_{[2n\times E_{Z}\sup _{g\in\cG'}\vert V_{n,Z}(g)\vert,\infty)}(\frac{\sqrt{n}d_n}{16}),
\end{align*}
or equivalently,
\begin{align*}
\log N\Big(\frac{d_n}{16},\mathcal{G}',\Vert .\Vert_n\Big)\leq &K\Big(\frac{32nE_{Z}\sup _{g\in\cG'}\vert V_{n,Z}(g)\vert}{\sqrt{n}d_n}\Big)^{\alpha}\mathbf{1}_{(0,32\sqrt{n} E_{Z}\sup _{g\in\cG'}\vert V_{n,Z}(g)\vert]}(d_n)+\\
&K\Big(\frac{32nE_{Z}\sup _{g\in\cG'}\vert V_{n,Z}(g)\vert}{\sqrt{n}d_n}\Big)^{2}\mathbf{1}_{[32\sqrt{n}E_{Z}\sup _{g\in\cG'}\vert V_{n,Z}(g)\vert,\infty)}(d_n).
\end{align*} 
Take $\varepsilon_i=Z_i/\sigma_\alpha=h(Z_i)$ for $i=1,...,n$, then $\mbox{var}(\varepsilon_i)=1$ and,
\begin{align*}
E_{\varepsilon}(\varepsilon_i)=E_{\varepsilon}(h(Z_i))=\int h(Z_i)\pi_\alpha(Z_i)dZ_i
=\frac{1}{\sigma_\alpha}\int Z_i\pi_\alpha(Z_i)dZ_i
=\frac{1}{\sigma_\alpha}E_Z(Z_i).
\end{align*} 
Therefore, $E_{Z}\sup _{g\in\cG'}\vert V_{n,Z}(g)\vert=\sigma_\alpha E_{\varepsilon}\sup _{g\in\cG'}\vert V_{n,\varepsilon}(g)\vert$ and,
\begin{align*}
\log N\Big(\frac{d_n}{16},\mathcal{G}',\Vert .\Vert_n\Big)\leq &K\Big(\frac{32n\sigma_\alpha  E_{\varepsilon}\sup _{g\in\cG'}\vert V_{n,\varepsilon}(g)\vert}{\sqrt{n}d_n}\Big)^{\alpha}\mathbf{1}_{(0,32\sigma_\alpha\sqrt{n} E_{\varepsilon}\sup _{g\in\cG'}\vert V_{n,\varepsilon}(g)\vert]}(d_n)+\\
&K\Big(\frac{32n\sigma_\alpha E_{\varepsilon}\sup _{g\in\cG'}\vert V_{n,\varepsilon}(g)\vert}{\sqrt{n}d_n}\Big)^{2}\mathbf{1}_{[32\sigma_\alpha\sqrt{n} E_{\varepsilon}\sup _{g\in\cG'}\vert V_{n,\varepsilon}(g)\vert,\infty)}(d_n),\\
\leq &K(32\sigma_\alpha)^\alpha n^{\frac{\alpha}{2}}\Big(\frac{\Ee\sup _{g\in\cG'}\vert V_{n,\varepsilon}(g)\vert}{d_n}\Big)^{\alpha}\mathbf{1}_{(0,32\sigma_\alpha\sqrt{n} E_{\varepsilon}\sup _{g\in\cG'}\vert V_{n,\varepsilon}(g)\vert]}(d_n)+\\
&K\Big(\frac{32\sigma_\alpha\sqrt{n}\Ee\sup _{g\in\cG'}\vert V_{n,\varepsilon}(g)\vert}{d_n}\Big)^{2}\mathbf{1}_{[32\sigma_\alpha\sqrt{n} E_{\varepsilon}\sup _{g\in\cG'}\vert V_{n,\varepsilon}(g)\vert,\infty)}(d_n),\\
\leq &K(32\sigma_\alpha)^\alpha n^{\frac{\alpha}{2}}\Big(\frac{\Ee\sup _{g\in\cG'}\vert V_{n,\varepsilon}(g)\vert}{d_n}\Big)^{\alpha}\mathbf{1}_{(0,32\sigma_\alpha\sqrt{n} E_{\varepsilon}\sup _{g\in\cG'}\vert V_{n,\varepsilon}(g)\vert]}(d_n)+\\
&K\mathbf{1}_{[32\sigma_\alpha\sqrt{n} E_{\varepsilon}\sup _{g\in\cG'}\vert V_{n,\varepsilon}(g)\vert,\infty)}(d_n).
\end{align*} 

\noindent
\underline{Proof of Step 4:}\label{ProofS4}

This Step consists in bounding from above the quantity $\mathbb{E}_\varepsilon  \sup_{g\in \mathcal{G}'} 
\vert V_{n,\varepsilon}(g)\vert$. According to Inequality \eref{but} we have,
\begin{align*}
\sum_{v\in \mathcal{P}}\vert V_{n,\varepsilon}(g_v)\vert \leq 
\kappa\Big\{ \sum_{v\in \mathcal{P}} \lambda_{n,v}^2 \| g_v\|_{\mathcal{H}_v}   + \sum_{v\in \mathcal{P}} \lambda_{n,v}  \| g_v\|_n \Big\},
\end{align*}
with $\lambda_{n,v}$ defined by Equation~\eref{lambda.eq} satisfying Equation \eref{condmu} for all $v\in \cP$. It follows
\begin{align*}
\sup_{g\in \mathcal{G}'} \sum_{v \in \cP}\vert V_{n,\varepsilon}(g_v)\vert 
&\leq 
\kappa\sup_{g\in \mathcal{G}'}\Big\{ \sum_{v\in \mathcal{P}}\lambda_{n,v}^2 \| g_v\|_{\mathcal{H}_v}  +\sum_{v\in \mathcal{P}} \lambda_{n,v}  \| g_v\|_n  \Big\},\\
&\leq  \frac{\kappa}{C_{1}}\sup_{g\in \mathcal{G}'} \Big\{ \sum_{v\in \mathcal{P}}\mu_v\| g_v\|_{\mathcal{H}_v}+\sum_{v\in \mathcal{P}} \gamma_{v} \| g_v\|_n\Big\}.
\end{align*}

Thanks to Condition \textbf{C1} and using $\| g_v\|_{\mathcal{H}_v} \leq 2$ we obtain then:
\begin{align*}
\sup_{g\in \mathcal{G}'} \sum_{v \in \cP}|V_{n,\varepsilon}(g_v)|
 &\leq
\frac{4\kappa}{C_{1}}\Big\{\sup_{g\in \mathcal{G}'} \sum_{v \in{S_f}}\mu_{v} \|g_{v}\|_{{\mathcal H}_{v}}
+  \sup_{g\in \mathcal{G}'} \sum_{v \in{S_f}}\gamma_{v}\| g_v\|_n\Big\},\\
 &\leq
\frac{4\kappa}{C_{1}}
\Big\{2 \sum_{v \in{S_f}}\mu_{v} + \sup_{g\in \mathcal{G}'} \sum_{v \in{S_f}}\gamma_{v}\| g_v\|_n\Big\}.
\end{align*}
Now, according to Condition \textbf{C2}, we get
\begin{align*}
\sup_{g\in \mathcal{G}'} \sum_{v \in \cP}\vert V_{n,\varepsilon}(g_v)\vert 
&\leq   
\frac{4\kappa}{C_1} \Big\{ 2\sum_{v\in {S_f}} \mu_{v} + 2\sup_{g\in \mathcal{G}'}\sum_{v\in {S_f}} \gamma_{v}\| g_v\|_2 +\sum_{v\in {S_f}} \gamma^2_{v}\Big\},\\
 &\leq  
 \frac{4\kappa}{C_1}  \Big\{\sum_{v\in {S_f}} (2\mu_v+\gamma^2_{v})+ 2\sup_{g\in \mathcal{G}'} (\sum_{v\in {S_f}}\gamma^2_{v})^{1/2} (\sum_{v\in {S_f}}\| g_v\|_2^2)^{1/2} \Big\},\\
 &\leq 
\frac{4\kappa}{C_1}\Big\{\sum_{v\in {S_f}} (2\mu_v+\gamma^2_{v})+2(\sum_{v\in {S_f}}\gamma^2_{v})^{1/2}d_n\Big\},
\end{align*}
where in the second inequality we used Cauchy Schwarz inequality and the third inequality coming from the fact that for all $g \in \cG'$,
$\|g\|^{2}_{2} = d_n^{2} \geq \sum_{v\in {S_f}} \| g_v\|_2^2 $. 


\hfill $\Box$

\subsection{Proofs of intermediate Lemmas\label{ProofsIntLemm.st}}
\subsubsection{Proof of Lemma  \ref{lemcomplex}\label{Prooflemcomplex}}
The kernel $k_v$ is written as :
\begin{equation*}
 k_{v} (X_{v}, X'_{v}) = \sum_{\ell\geq 1} \omega_{v,\ell} 
\phi_{v,\ell}(X_{v})\phi_{v,\ell}(X'_{v})
\end{equation*}
where $\{\phi_{v,\ell}\}_{\ell=1}^{\infty}$
is an orthonormal basis of $L^2(P_{v})$ with $P_v
= \prod_{a \in v} P_{a}$.

Let us consider the class of functions $\cK(t)$ defined as 
\begin{equation*}
  \cK(t) =\left\{ g_v \in \cH_{v}, \|g_v\|_{\cH_v}\leq 2, \|g_v\|_{2}\leq t \right\}.
\end{equation*}
It comes that
\begin{equation*}
g_{v} = \sum_{\ell} a_{\ell} \phi_{v,\ell}, 
\:\: \mbox{ with }
\|g_{v}\|^{2}_{\cH_{v}} = \sum_{\ell} \frac{a_{\ell}^{2}}{\omega_{v,\ell} }
\leq 4, \mbox{ and }
\|g_{v}\|^{2}_{2} = \sum_{\ell} a_{\ell}^{2} \leq t^{2}
\end{equation*}

In the following, we set  $\mu_{v,\ell} (t)= \min\left\{t^2, \omega_{v,\ell}\right\}$.  Hence
\begin{equation}
 \sum_{\ell} \frac{a_{\ell}^{2}}{\mu_{v,\ell}(t)} 
\leq 
  \frac{1}{t^{2}} \sum_{\ell}a_{\ell}^{2} + \sum_{\ell} \frac{a_{\ell}^{2}}{\omega_{v,\ell}}
= \frac{1}{t^{2}} \|g_{v}\|_2^{2} + 
\|g_{v}\|_{{\mathcal H}_{v}}^{2} \leq 5,
\label{ineq36.eq}
\end{equation}
as soon as  $g_{v}\in \cK(t)$.

Now, let us prove the lemma:
\begin{align*}
\EXe W_{n,2,v}(t)&=
\EXe \sup_{g \in \cK(t) }
\vert \frac{1}{n}
\sum_{i=1}^n\varepsilon_i \sum_{\ell}a_{\ell}\phi_{v,\ell}(X_{vi}) 
\vert, \\
&= \EXe  \sup_{g \in \cK(t)}
\vert 
\frac{1}{n}\sum_{\ell}\frac{a_{\ell}}{\sqrt{\mu_{v,\ell}(t)}} \sum_{i=1}^n\varepsilon_i
  \sqrt{\mu_{v,\ell} (t)}\phi_{v,\ell}(X_{vi})\vert, \\
& \leq  \sqrt{5}\sqrt{
\EXe  \sum_{\ell} \Big( \frac{1}{n}\sum_{i=1}^n\varepsilon_i
  \sqrt{\mu_{v,\ell}(t) }\phi_{v,\ell}(X_{vi})\Big)^{2}} .
  \end{align*}
The last inequality follows from the Cauchy-Schwartz inequality and Inequality~\eref{ineq36.eq}.
  Now, simple calculation leads to 
\begin{equation*}
 \EXe  W_{n,2,v}(t) \leq \sqrt{5} \sqrt{\frac{1}{n}\sum_{\ell} \mu_{v,\ell} (t) }.
\end{equation*}

\hfill $\Box$

\subsubsection{Proof of Lemma \ref{lemcompnormes1}\label{Prooflemcompnormes1}}
Using that $\vert \sqrt{a}-\sqrt{b}  \vert  \leq \sqrt{\vert a-b\vert},$
we get
$$\left\vert\| g_v\|_2- \| g_v\|_n\right\vert   \leq \sqrt{\left\vert \| g_v\|_2^2- \| g_v\|_n^2\right\vert} .$$
Hence
$$\left\lbrace  \| g_v\|_\infty\leq b , \;
\left\vert\| g_v\|_2- \| g_v\|_n\right\vert \geq \frac{b t}{2}\right\rbrace
\subset \left\lbrace \left\vert\| g_v\|_2^2- \| g_v\|_n^2\right\vert \geq \frac{b^2 t^2}{4}\right\rbrace.$$
The centered process
\begin{align*}
\left\vert\| g_v\|_2^2- \| g_v\|_n^2\right\vert=\vert\frac1n \sum_{i=1}^n g_v^2(X_{v,i})-\mathbb{E}(g_v^2(X_v))\vert,
\end{align*}
satisfies a concentration inequality given, for example,  by  Theorem
2.1 in \cite{bartlett2005} : if $\mathcal{C}$ is a class of
functions $f$ such that $\|f\|_{\infty} \leq B$ and $E f(X)=0$, and
if there exists $\gamma>0$ 
such that for every $f \in \mathcal{C}$, $\text{Var} f(X) \leq \gamma^{2}$. Then for every
$x>0$, with probability at least $1-e^{-x}$, 
\begin{equation}
 \sup_{f\in \mathcal{C}}\frac{1}{n}\vert \sum_{j=1}^n
  f(X_j)\vert \leq \inf_{\alpha>0}\Big\{
2(1+\alpha) E (\sup_{f\in \mathcal{C}}\frac{1}{n}\vert \sum_{j=1}^n
  f(X_j)\vert)
+ \sqrt{\frac{2  x}{n}} \gamma + B \Big(\frac{1}{3} +
  \frac{1}{\alpha}\Big)\frac{x}{n}
\Big\}.
\label{Bartlet.eq}
\end{equation}

For any $t>0$, for $\cG(t)$ defined by~\eref{calG.eq},
let us consider  the class of functions $\cC(t)$
defined as follows
\begin{align*}
 \mathcal{C}(t) =\Big\{f \mbox{ such that } f=g_v^2-\mathbb{E}(g_v^2), \mbox{
  with } g_v\in \cG(t) \Big\}. 
\end{align*}
  Note that if $f \in \mathcal{C}(t)$, $\EX f(X_v)=0$ and $\|f\|_{\infty}
\leq  b^{2}$. We have to study 
\begin{align*}
 \gamma^{2}(t)  = \sup_{g_v \in \cG(t)} 
\EX \left(   g_v^2(X)-\|g_v\|^{2}_{2})\right)^{2} \mbox{ and }
\Gamma(t)  =  \EX \Big( \sup_{g_v \in \cG(t)} \left|
    \|g_v\|^{2}_{n}-\|g_v\|^{2}_{2}\right| \Big).
\end{align*}
It is easy to see that 
\begin{align*}
 \gamma^{2} (t)
 \leq   b^{2} \sup_{g_v \in \cG(t)} \EX \left(
  g_v(X)+\|g_v\|_{2}\right)^{2} \leq 4 b^{2}t^{2}.
\end{align*}

Let  $\zeta_i$ be i.i.d. Rademacher random variables and let $E_{X,\zeta}$ denotes the
expectation with respect to the law 
of $(X,\zeta)$.
By  a symmetrization argument, 
$$\Gamma(t)  \leq 2 E_{X,\zeta}
\sup_{g_v\in \mathcal{G}(t)}
\vert\frac{1}{n} \sum_{i=1}^n \zeta_i g_v^2(X_i)\vert .$$

Since $\| g_v\|_{\infty}\leq b$, applying the  contraction principal (see \cite{LedouxTal:91}) we get that, for $Q_{n,v}(t)$ defined by
\eref{Qn},
\begin{align*}
E_{X,\zeta}
\sup_{g_v\in \mathcal{G}(t)}
\vert\frac{1}{n} \sum_{i=1}^n \zeta_i g_v^2(X_i)\vert
\leq 4 b E_{X,\zeta}
\sup_{g_v\in \mathcal{G}(t)}
\vert\frac{1}{n} \sum_{i=1}^n \zeta_i g_v(X_i)\vert
\leq  4 b Q_{n,v}(t).
\end{align*}
The last inequality was 
proved by \cite{Mendelson2002}, Theorem 41 (see the proof of
Lemma~\ref{lemcomplex}).  Now, thanks
to~\eref{Bartlet.eq} we get that 
for all $x>0$, with probability greater than $1-e^{-x}$
\begin{equation*}
 \sup_{g_v \in \cG(t)}\vert \|g_{v}\|_{n}^{2} - \|g_{v}\|_{2}^{2}
  \vert \leq \inf_{\alpha>0}\Big\{
16 (1+\alpha)  b Q_{n,v}(t)
+ \sqrt{\frac{2  x}{n}}2 b t + b^{2} \left(\frac{1}{3} +
  \frac{1}{\alpha}\right)\frac{x}{n}
\Big\}.
\end{equation*}
Taking $x=c_{2} n t^{2}$, $t \geq  \nu_{n}$, we have that with
probability greater than $1-e^{-c_{2} n t^{2}}$
\begin{equation*}
 \sup_{g_v \in \cG(t)}\left\vert \|g_{v}\|_{n}^{2} - \|g_{v}\|_{2}^{2}
  \right\vert \leq \inf_{\alpha>0}t^{2} \Big\{
16(1+\alpha) b \Delta 
+ \sqrt{2  c_{2}}4 b  + b^{2} \left(\frac{1}{3} +
  \frac{1}{\alpha}\right)c_{2}  
\Big\}.
\end{equation*}
The infimum of the right hand side is reached in $\alpha=\sqrt{c_{2}
  b/ 16 \Delta}$,
and equals  
\begin{equation*}
 \frac{b^{2} c_{2}}{3} + 8 \sqrt{\Delta c_{2}} b^{3/2} + 4(4 \Delta + 
 \sqrt{2  c_{2}}) b.
\end{equation*}
The constants $\Delta$ and $c_{2}$ should satisfy that this infimum is
strictly smaller than $b^{2}/4$. 
For example, if $16 \Delta <b/8$, it remains to choose
$c_{2}$ small enough such that 
\begin{equation*}
b\left(\frac{c_{2}}{3} + \frac{\sqrt{2 c_{2}}}{2}\right) + 4\sqrt{2 c_{2}}
 < \frac{b}{8}.
\end{equation*}

\hfill \textbf{$\Box$}

\subsubsection{Proof of Lemma \ref{lemcompnormes2}\label{Prooflemcompnormes2}}
Let $t> \nu_{n,v}$ and $h$ be defined as $$h =\frac{ t g_v }{ \|g_v\|_{2}}.$$ If
$g_v$ satisfies the assumptions of the lemma, then $h$  satisfies $\| h\|_2=t$, $\|h\|_{\mathcal{H}}\leq 2$ and $\|h\|_{\infty} \leq b$.
Applying Lemma \ref{lemcompnormes1}
(page~\pageref{lemcompnormes1}) to the function  $h$, we
obtain that for all $ t \geq \nu_{n,v}$, with probability
greater than $1 - \exp(-c_2 n t^2)$, we have
$$ |t - \|h\|_{n}| \leq \frac{bt}{2}\quad\text{for all}\quad h \in \cG(t).$$ This concludes the proof of the lemma.

\hfill \textbf{$\Box$}

\subsubsection{Proof of Lemma \ref{concentration1}}
\label{Proofconcentration1}
We apply Corollary \ref{lipschitz} to $$\phi(\varepsilon_{1}, \ldots,
\varepsilon_{n}) =\frac{\sqrt{n}}{t} W_{n,n,v}(t).$$ Using Cauchy-Schwarz Inequality and the fact that $\Vert g_v\Vert_n\leq t$,
\begin{align*}
| \phi(\varepsilon) - \phi(\varepsilon')| \leq \frac{\sqrt{n}}{t}\sup_{\|g_v\|_{n} \leq t}\|g_v\|_{n} \|\varepsilon-\varepsilon' \|_{n}
\leq   \frac{\sqrt{n}}{t}t \|\varepsilon-\varepsilon' \|_{n},
\end{align*} 
leading to $\|\phi\|_{L}=1$. So,
\begin{align*}
\PXe \Big(\vert \frac{\sqrt{n}}{t}W_{n,n,v}( t)-\frac{\sqrt{n}}{t}E_{\varepsilon} W_{n,n,v}( t)\vert \geq u\Big) \leq 2B \exp\Big(-\frac{u^2}{8A}\Big),
 \end{align*}
and Lemma
  \ref{concentration1} is proved by taking $\delta=u/\sqrt{n}$.
  
\hfill $\Box$

\subsubsection{Proof of Lemma \ref{concentration2}}
\label{Proofconcentration2} 
We start with the proof of \eref{concentn2} in Lemma
\ref{concentration2} by applying once again Corollary \ref{lipschitz}, to
the function $$\phi(\varepsilon) = \phi(\varepsilon_{1}, \ldots,
\varepsilon_{n}) = \frac{\sqrt{n}}{2t} W_{n,2,v}(t).$$
On the event  $\Omega_{v, t}$ defined by~\eref{Omega},
we have $$\|g_{v}\|_{n} \leq \frac{bt}{2} + \|g_{v}\|_{2}.$$ Besides if
$\|g_{v}\|_{\cH_{v}} \leq 2$, then $\|g_{v}\|_{\infty} \leq
2$. Therefore applying Lemma~\ref{lemcompnormes1} with $b=2$, we get
that if $\|g_{v}\|_{2} \leq t$,
\begin{align*}
| \phi(\varepsilon) - \phi(\varepsilon')| \leq \frac{\sqrt{n}}{2t}\sup_{\|g_v\|_{n} \leq 2t}\|g_v\|_{n} \|\varepsilon-\varepsilon' \|_{n}
\leq  \frac{\sqrt{n}}{2t}2 t \|\varepsilon-\varepsilon' \|_{n},
\end{align*} 
leading to $\|\phi\|_{L}=1$. So,
$$\PXe \Big(\Big\lbrace\vert \frac{\sqrt{n}}{2t}W_{n,2,v}( t)-
\frac{\sqrt{n}}{2t}E_\varepsilon (W_{n,2,v}( t))
\vert    \geq u\Big\rbrace \cap \Omega_{v,t}^c\Big) \leq 2B \exp\Big(-\frac{u^2}{8A}\Big),$$
and inequality \eref{concentn2} in Lemma  \ref{concentration2} is proved by taking $\delta=2u/\sqrt{n}$.

\medskip
\noindent
We now come to the proof of the inequality \eref{concentration} in Lemma
\ref{concentration2} 
using a Poissonian inequality for self-bounded processes (see \cite{Boucheron2000ASC}) and  Theorem 5.6, p 158 in \cite{massart2007concentration}). Let us recall it in the particular case we are interested in:
\begin{theo}
\label{autoborne}
Let  $X_1,\cdots, X_n$  be $n$ i.i.d. random variables. 
For
$i\in \{1,\cdots,n\}$ let  $$X_{(-i)}=(X_1,\ldots, X_{i-1},X_{i+1},\ldots,X_n).$$
Let $h$ be a non-negative and bounded measurable function of $X=(X_1,\cdots,X_n)$. Assume that for all $i\in\{1,\cdots,n\}$, there exists a measurable function
$h_i$ of  $X_{ (-i)}$ such that
$0<h-h_i\leq 1,$ and 
 $\sum_{i=1}^n (h-h_i)\leq h$.
 Then, for 
all  $x>0$, we have $$P \Big(h\geq E(h) +x \Big)  \leq \exp\Big(
  -\frac{x^2}{2 E(h)}   \Big).$$
\end{theo}
We apply this result to $h$ defined as 
\begin{equation*}
h=h(X_1,\cdots,X_n)= 
n \Ee W_{n,2,v}(t) = 
n \Ee
\sup \Big\{\vert V_{n,\varepsilon}(g_v)\vert,\: \| g_v\|_2\leq t,\: 
\| g_v\|_{\cH_v}\leq 2\Big\}   .
\end{equation*}
The variable $h$ is positive, and because the distribution of $(\varepsilon_{1}, \ldots,\varepsilon_{n})$ is symmetric, we have that 
\begin{equation*}
h=E_{\varepsilon} 
\sup \Big\{ n V_{n,\varepsilon}(g_v),\: \| g_v\|_2\leq t, \:
\| g_v\|_{\cH_v}\leq 2\Big\}.
\end{equation*}
Let $\tau$ be the function in $\mathcal{H}_v$ such that 
$h= E_{\varepsilon} n V_{n,\varepsilon}(\tau)$ (note that $\tau$
depends on $(X_1, \ldots, X_n)$ and on $(\varepsilon_{1}, \ldots,
\varepsilon_{n})$), and let 
\begin{equation*}
h_{i}=  \Ee  \sup_{g_v} 
\sum_{j\not =i} \varepsilon_j g_{v}(X_j).
\end{equation*}
 We  show that $h$ and $h_{i}$ satisfy the assumptions of Theorem~\ref{autoborne}: 
\begin{align*}
 h-h_{i} &=  E_{\varepsilon} \Big(\varepsilon_i \tau( X_i) +
\sum_{j \neq i}\varepsilon_j \tau( X_j)
- \sup_{g_{v}} \sum_{j \neq i} \varepsilon_jg_{v}( X_j) \Big), \\
&\leq   E_{\varepsilon}  \Big(\varepsilon_i \tau( X_i)\Big),\\ 
&\leq  E_{\varepsilon}  \Big(\vert\varepsilon_i\vert \sup_{x \in \cX}|\tau(X)|\Big),\\
&\leq  2E_{\varepsilon} \Big( \vert\varepsilon_i\vert\Big),
\end{align*}
where the last inequality comes from the fact that $\sup_{x \in
  \cX}|\tau(X)| \leq \|\tau\|_{\cH_v} \leq 2$. 

Let $Z=(Z_1,...,Z_n)$ be i.i.d. random variables distributed with density $\pi_\alpha\in\mathcal{D}$ defined in Equation \eref{familypi} with $\sqrt{\mbox{var}(Z_i)}=\sigma_\alpha$. Take $\varepsilon_i=Z_i/\sigma_\alpha$ for $i=1,...,n$, then $\mbox{var}(\varepsilon_i)=1$ and  $ E_{\varepsilon}(\vert\varepsilon_i\vert)=E_Z(\vert Z_i\vert)/\sigma_\alpha$. We have:
\begin{align*}
E_{Z}(\vert Z_i\vert)=\int_\mathbb{R}\vert Z_i\vert a_\alpha \exp\Big({-\vert Z_i\vert}^\alpha\Big)dZ_i.
\end{align*}
Take $\vert Z_i\vert=u^{1/\alpha}$ and 
\begin{align*}
dZ_i=\left\{ \begin{array}{rcl}\frac{1}{\alpha}u^{\frac{1}{\alpha}-1}du & \mbox{if} &Z_i\geq 0,\\
-\frac{1}{\alpha}u^{\frac{1}{\alpha}-1}du& \mbox{if} &Z_i\leq 0.
\end{array}\right.
\end{align*}
Therefore,
\begin{align*}
E_{Z}(\vert Z_i\vert)&=\int_0^{+\infty} a_\alpha u^{\frac{1}{\alpha}}\exp(-u)\frac{1}{\alpha}u^{\frac{1}{\alpha}-1}du-\int_{+\infty}^0 a_\alpha u^{\frac{1}{\alpha}}\exp(-u)\frac{1}{\alpha}u^{\frac{1}{\alpha}-1}du,\\
& =2\int_0^{+\infty} a_\alpha u^{\frac{1}{\alpha}}\exp(-u)\frac{1}{\alpha}u^{\frac{1}{\alpha}-1}du,\\
& =a_\alpha\int_0^{+\infty} \frac{2}{\alpha}u^{\frac{2}{\alpha}-1}\exp(-u)du,\\
&=a_\alpha \frac{2}{\alpha}\Gamma(\frac{2}{\alpha})=a_\alpha\Gamma(1+\frac{2}{\alpha}),
\end{align*}
where $\Gamma(.)$ is the gamma function.

It follows that,
\begin{equation*}
h-h_{i}\leq \frac{2a_\alpha}{\sigma_\alpha}\Gamma(1+\frac{2}{\alpha}).
\end{equation*}
Moreover, $h-h_i \geq 0$ since 
\begin{align*}
h  =  E_\varepsilon \Big(\sup_{g_{v}} \sum_{j=1}^{n} \varepsilon_j g_v( X_j)\Big)
 =  E_\varepsilon  \Big(E_{\varepsilon_{i}} \sup_{g_{v}}
\sum_{j=1}^{n} \varepsilon_j g_v( X_j) \Big)
\geq   E_\varepsilon  \Big(\sup_{g_{v}}
E_{\varepsilon_{i}}\sum_{j=1}^{n} \varepsilon_j g_v( X_j) \Big)= h_{i}.
\end{align*}
Finally we have: 
$$\sum_i (h-h_{i})= \sum_{i=1}^n E_\varepsilon   
\Big( \varepsilon_i\tau(X_i)+ \sum_{j\not=i}^n \varepsilon_j
  \tau(X_j) 
-\sup_{g-v}\sum_{j\not =i}^n \varepsilon_j g_v(X_j) \Big)
\leq \sum_{i=1}^n \Ee  \varepsilon_i\tau(X_i) = h.$$
Therefore, following Theorem~\ref{autoborne}, we get that for all
postive $u$
\begin{align*}
P_{X,\varepsilon} \Big( E_\varepsilon  W_{n,2,v}(t)
    -E_{X,\varepsilon} W_{n,2,v}(t) \leq  \frac{u}{n} \Big)\leq \exp\Big( -\frac{u^2}{
E_{X,\varepsilon} W_{n,2,v}(t)} \Big).
\end{align*}
As $E_{X,\varepsilon}  W_{n,2,v}(t)\leq Q_{n,v}(t)$, see Lemma~\ref{lemcomplex} page~\pageref{lemcomplex}, we get the expected result since for all
positive $x$
\begin{align*}
P_X \Big(E_\varepsilon W_{n,2,v}(t)\geq E_{X,\varepsilon} W_{n,2,v}(t)+x\Big)\leq \exp\Big(-\frac{n x^2}{ Q_{n,v}(t)} \Big).
\end{align*}

\hfill $\Box$

\subsubsection{Proof of Lemma \ref{Case1}}
\label{ProofCase1}  
From Lemma \ref{concentration1}, page~\pageref{concentration1} with $t=\lambda_{n,v}=\delta$, with  probability greater than $1-2B\exp(-n\lambda_{n,v}^2/8A)$, we get that:
 \begin{equation}
 E_\varepsilon(W_{n,n,v}(\lambda_{n,v}))\leq \lambda^2_{n,v}+W_{n,n,v}(\lambda_{n,v})
 \end{equation}
 
The next step consists in comparing $W_{n,n,v}( \lambda_{n,v})$ and $
W_{n,2,v}( 2\lambda_{n,v})$. Recall that $\lambda_{n,v} \geq \nu_{n,v}$, see~\eref{lambda.eq}. Let $g_{v}$ such that $\|g_v\|_{n}\leq \lambda_{n,v}$. 
\begin{itemize}
\item When $\| g_v\|_2 \leq \lambda_{n,v}$, according to Lemma \ref{lemcompnormes1} (page~\pageref{lemcompnormes1}), taking $b=2$ , since since $\| g_v\|_n\leq \lambda_{n,v}$, we get that  with probability greater than $1-\exp(-c_{2} n \lambda_{n,v}^{2})$, 
$$\| g_v\|_n -\lambda_{n,v}\leq \| g_v\|_2 \leq
\|g_v\|_n+\lambda_{n,v} \leq 2 \lambda_{n,v}.$$
\item When $\| g_v\|_2 \geq t$, we apply Lemma \ref{lemcompnormes2} (page~\pageref{lemcompnormes2}) with $b=2$.
For any function $g_{v}$ such that $\| g_{v}\|_\infty \leq 2$, and $\| g_{v}\|_2\geq \lambda_{n,v}$, we have $\| g_{v}\|_2\leq 2 \| g_{v}\|_n \leq 2 \lambda_{n,v}$.
\end{itemize}
This implies that,  with probability greater than $1-\exp(-c_2n \lambda_{n,v}^2)$
 we have
$$W_{n,n,v}( \lambda_{n,v}) \leq W_{n,2,v}( 2\lambda_{n,v}).$$

We now study the process $W_{n,2,v}( \lambda_{n,v})$. By applying
\eref{concentn2} in Lemma  \ref{concentration2},
page~\pageref{concentration2}, with $\delta=t=\lambda_{n,v}$ we get
that with  probability greater than $1-2B\exp(- n\lambda_{n,v}^2/32A)$
$$ W_{n,2,v}( \lambda_{n,v})\leq \lambda_{n,v}^2+ E_\varepsilon (W_{n,2,v}( \lambda_{n,v})).$$

It follows that
\begin{align*}
\Ee W_{n,n,v}( \lambda_{n,v})&\leq \lambda_{n,v}^2+W_{n,n,v}( \lambda_{n,v} )),\\
&\leq \lambda_{n,v}^2+W_{n,2,v}( 2\lambda_{n,v})),\\
&\leq 5\lambda_{n,v}^2+\Ee  (W_{n,2,v}( 2\lambda_{n,v})).
\end{align*}
Next, we apply \eref{concentration} in Lemma \ref{concentration2}, with
$t=2\lambda_{n,v}$ and $x= 4 \lambda^{2}_{n,v}$. We get that $$
\Ee W_{n,2,v}( 2\lambda_{n,v})\leq
4\lambda_{n,v}^2+ \EXe  (W_{n,2,v}(
2\lambda_{n,v})),$$
with probability greater than 
\begin{equation*}
 1 -2 \exp(-16 \frac{n \lambda^{4}_{n,v}}{ Q_{n,v}(2
     \lambda_{n,v})})
\geq 1 -2 \exp(- \frac{4 n \lambda^{2}_{n,v}}{\Delta }).
\end{equation*}
The last inequality comes from the definition of $\nu_{n,v}$,
see~\eref{nu}, and from the fact that $\lambda_{n,v} \geq \nu_{n,v}$,
see~\eref{lambda.eq}.

Putting everything together, we get that with probability greater than
$1 - c_{1}\exp(-c_2 n\lambda_{n,v}^2)$ for some positive constants
$c_1, c_2$, 
\begin{align*}
\Ee W_{n,n,v}( \lambda_{n,v})
&\leq 9\lambda_{n,v}^2+ \EXe  (W_{n,2,v}( 2\lambda_{n,v})),\\
&\leq 9\lambda_{n,v}^2+  Q_{n,v}(2\lambda_{n,v}),
\mbox{ thanks to Lemma~\ref{lemcomplex}, page~\pageref{lemcomplex}},
\\
&\leq 9\lambda_{n,v}^2+ 4 \Delta \lambda_{n,v}^2.
\end{align*}
Applying once again  Lemma \ref{concentration1},
page~\pageref{concentration1}, we get that
\begin{equation*}
W_{n,n,v}(\lambda_{n,v}) \leq   \Ee  W_{n,n,v}(
\lambda_{n,v})+ \lambda_{n,v}^2 \leq   \Big(10+4 \Delta\Big)\lambda_{n,v}^2. 
\end{equation*}
This ends the proof of the lemma by taking $\kappa=10+4 \Delta$.

\hfill $\Box$

\section{Proof of Corollary \ref{oracle2}}\label{prooforacle2}
According to Theorem \ref{oracle} we have with high probability,
\begin{equation}
\label{resoracle}
 \|  \widehat{f}-m\|^{2}_{n} \leq
C \inf_{f \in \cF}\Big\{  \| m - f\|^{2}_{n}
+  \sum_{v \in S_{f}} (\mu_{v} + \gamma_{v}^{2})
 \Big\}.
\end{equation}
Besides, for all $K>0$,
\begin{align}
\label{fmdec}
\Vert \widehat{f}-m\Vert_2^2\leq(1+K)\Vert \widehat{f}-f\Vert_2^2+(1+\frac{1}{K})\Vert m-f\Vert_2^2.
\end{align}
We consider once again two cases defined in page \pageref{3Cases}.

\underline{Case 1:} $\Vert \widehat{f}-f\Vert_2\leq\Vert \widehat{f}-f\Vert_n$,

In this case Equation \eref{fmdec} gives,
\begin{align*}
\Vert \widehat{f}-m\Vert_2^2\leq(1+K)\Vert \widehat{f}-f\Vert_n^2+(1+\frac{1}{K})\Vert m-f\Vert_2^2.
\end{align*}
Then, using Equations \eref{trick2} and \eref{resoracle} we obtain the result.

\underline{Case 2:} $\Vert \widehat{f}-f\Vert_2\geq\Vert \widehat{f}-f\Vert_n$,

Apply  Lemma \ref{lemcompnormes3} (page~\pageref{lemcompnormes3}) and conclude that conditioning on the events $\mathcal{T}$ and $\mathcal{A}$, defined by \eref{evtTau} and \eref{A},
 then $\widehat{f}-{f}$ belongs to $\mathcal{G}(f)$ defined in Lemma \ref{lemcompnormes3}. 
Now, conditioning on the event $\cC$
we get the result as in Case 1 since, 
$$\|  \widehat{f}-{f}\|_2\leq \sqrt{2}\| \widehat{f}-{f}\| _n.$$

\hfill $\Box$

\section{Proofs of Section \ref{sec:sudakov}}\label{proofssudakov}
\subsection{Proof of Remark \ref{ubrelation}}\label{proofubrelation}
From Lemma \ref{lem23tala} we have $U_{\tilde{\alpha}}(u)\subset (u^{1/2}B_2 + u^{1/\tilde{\alpha}}B_{\tilde{\alpha}} )$. It suffices to show that $(u^{1/2}B_2 + u^{1/\tilde{\alpha}}B_{\tilde{\alpha}} )\subset 2\times \max (u^{1/2} ,u^{1/\tilde{\alpha}})B_2$. 

Consider $x\in u^{1/2}B_2+u^{1/\tilde{\alpha}}B_{\tilde{\alpha}}$, $x=y+z$ with $y\in u^{1/2}B_2$, means $\sum_{i=1}^ny_i^2\leq u$, and $z\in u^{1/\tilde{\alpha}}B_{\tilde{\alpha}}$, means $\sum_{i=1}^nz_i^{\tilde{\alpha}}\leq u$. Moreover, we know that $\Vert x\Vert\leq\Vert y\Vert+\Vert z\Vert$ which leads to $\Vert x\Vert\leq u^{1/2}+u^{1/\tilde{\alpha}}$ and $\Vert x\Vert^2\leq 2(u+u^{2/\tilde{\alpha}})\leq 4\times \max(u,u^{2/\tilde{\alpha}})$.

\hfill $\Box$

\subsection{Proof of Corollary \ref{corosudakov}}\label{proofcorosudakov}
From Equation (\ref{sudkvmax}) we have $N(2\times \max(M^{1/2},M^{1/\tilde{\alpha}}),T,\Vert .\Vert)\leq \exp(KM)$. 

Using this on $sT$ for $s > 0$ we have $sM = E_{Z} \sup_{t'\in sT} \sum_{i=1}^nt'_iZ_i$ and, 
\begin{align*}
N(2\times \max((sM)^{1/2},(sM)^{1/\tilde{\alpha}}),sT,\Vert .\Vert)\leq \exp(KsM).
\end{align*}
Moreover, 
\begin{align*}
N(2\times \max((sM)^{1/2},(sM)^{1/\tilde{\alpha}}),sT,\Vert .\Vert)=N(\frac{2}{s}\times \max((sM)^{1/2},(sM)^{1/\tilde{\alpha}}),T,\Vert .\Vert),
\end{align*}
since for all $t_1,t_2\in T$ and some constant $C$, $\Vert st_1-st_2\Vert\leq C$ is equivalent to $\Vert t_1-t_2\Vert\leq C/s$. 

We obtain then,
\begin{align*}
N(\frac{2}{s}\times \max((sM)^{1/2},(sM)^{1/\tilde{\alpha}}),T,\Vert .\Vert)\leq \exp(KsM).
\end{align*}
As in Remark \ref{lemsudkvmax} for $u=sM$ we consider two following cases (recall that $1<\tilde{\alpha}< 2$):
\begin{itemize}
\item[(i)] If $sM\leq 1$ we have $(sM)^{1/\tilde{\alpha}}\leq (sM)^{1/2}$ and so,
\begin{align*}
N(2(\frac{M}{s})^{1/2},T,\Vert .\Vert)\leq \exp(KsM).
\end{align*}
Take $\delta =2(M/s)^{1/2}$ and thus $s=4M/\delta^2$. Moreover, $sM\leq 1$ (i.e. $(4M/\delta^2)\times M\leq 1$) and so $\delta\geq 2M$. Finally, we obtain in this case:
\begin{align*}
\forall \delta\geq 2M,\:\log N(\delta,T,\Vert .\Vert)\leq K(\frac{2M}{\delta})^2.
\end{align*}
\item[(ii)] If $sM\geq 1$ we have $(sM)^{1/2}\leq(sM)^{1/\tilde{\alpha}}$ and so,
\begin{align*}
N(\frac{2}{s}(sM)^{1/\tilde{\alpha}},T,\Vert .\Vert)\leq \exp(KsM).
\end{align*}
Take $\delta =(2/s)(sM)^{1/\tilde{\alpha}}$ and thus $s=(2/\delta)^{\tilde{\alpha}/(\tilde{\alpha}-1)}M^{1/(\tilde{\alpha}-1)}$. Moreover, $sM\geq 1$ (i.e. $(2M/\delta)^{\tilde{\alpha}/(\tilde{\alpha}-1)}\geq 1$) and so $0<\delta\leq 2M$. Finally, we obtain in this case:
\begin{align*}
\forall 0<\delta\leq 2M,\:\log N(\delta,T,\Vert .\Vert)\leq K(\frac{2M}{\delta})^{\tilde{\alpha}/(\tilde{\alpha}-1)}=K(\frac{2M}{\delta})^{\alpha}.
\end{align*}
\end{itemize}

\hfill $\Box$

\section{Proofs of Section \ref{sec:conineq}}\label{proofsec:conineq}
\subsection{Proof of Lemma \ref{lipsc2}}\label{prooflipsc2}  
In order to prove this Lemma it suffices to show that $\Pi_\alpha\in\mathcal{M}(m,\rho^2)$ for some $m$. To do so, we use Example \ref{exadamczak}.

First show $d\log \Pi_\alpha ([t,\infty))/dt\leq -t/\rho^2$:

We know that
\begin{align*}
\frac{d}{dt}\log \Pi_\alpha ([t,\infty))=\frac{d}{dt}\log(1-\Pi_\alpha((-\infty,t]))=-\frac{\pi_\alpha(t)}{1-\Pi_\alpha((-\infty,t])}=-\frac{\pi_\alpha(t)}{\Pi_\alpha([t,\infty))}.
\end{align*}
For all $t> 0$ we have,
\begin{align*}
\Pi_\alpha([t,\infty))=\int_t^\infty a_\alpha \exp(-\vert x\vert^\alpha)dx=\int_t^\infty a_\alpha \exp(- x^\alpha)dx.
\end{align*}
Take $x=u^{1/\alpha}$, so $dx= (1/\alpha)u^{(1/\alpha)-1} du$, and
\begin{align*}
\Pi_\alpha([t,\infty))&=\int_{t^\alpha}^\infty \frac{a_\alpha}{\alpha}u^{(1/\alpha)-1}\exp(-u)du,\\
&=\frac{a_\alpha}{\alpha}\Gamma(\frac{1}{\alpha},t^\alpha),
\end{align*}
where $\Gamma(\frac{1}{\alpha},t^\alpha)$ is incomplete gamma function. Moreover, for $s\in\mathbb{R}$ as $x\rightarrow \infty$,
\begin{align*}
\frac{\Gamma(s,x)}{x^{s-1}\exp(-x)}\rightarrow 1.
\end{align*}
 Therefore,
 $$\Pi_\alpha([t,\infty))=\frac{a_\alpha}{\alpha}t^{1-\alpha}\exp(-t^\alpha).$$
Since $t>0$ so $\pi_\alpha(t)=a_\alpha \exp(-t^\alpha)$, and 
\begin{align*}
\frac{d}{dt}\log \Pi_\alpha ([t,\infty))=-\frac{\alpha a_\alpha \exp(-t^\alpha)}{a_\alpha t^{1-\alpha} \exp(-t^\alpha)}=-\alpha t^{\alpha-1}.
\end{align*} 
The inequality $-\alpha t^{\alpha-1}\leq -t/\rho^2$ (i.e. $t^{\alpha-2}\geq 1/\alpha\rho^2$) holds for all $\alpha> 2$ and $t\geq (1/\alpha\rho^2)^{1/(\alpha-2)}$.

Second show $d\log \Pi_\alpha ((-\infty,-t])/dt\leq -t/\rho^2$:

The probability distribution $\Pi_\alpha$ is symmetric, therefore $\Pi_\alpha ((-\infty,-t])=\Pi_\alpha ([t,\infty))$, and
\begin{align*}
\frac{d}{dt}\log \Pi_\alpha ((-\infty,-t])=\frac{d}{dt}\log \Pi_\alpha ([t,\infty))=-\alpha t^{\alpha-1},
\end{align*}
which is smaller than $-t/\rho^2$ if $\alpha> 2$ and $t\geq (1/\alpha\rho^2)^{1/(\alpha-2)}$.

Take $m=(1/\alpha\rho^2)^{1/(\alpha-2)}$, then for $x\geq m$, $\Pi_\alpha$ verifies the Equations $(\ref{mfami})$. That is $\Pi_\alpha\in\mathcal{M}((1/\alpha\rho^2)^{1/(\alpha-2)},\rho^2)$. 

\hfill $\Box$

\subsection{Proof of Remark \ref{expecfini}}\label{proofexpecfini}
If $\alpha=2$, according to the Laplace transform of the Gaussian function we have 
\begin{align}
\label{glap}
E(\exp({s\vert Z\vert}))=2a_\alpha\sqrt{\pi} \exp\Big({\frac{s^2}{4}}\Big).
\end{align}
If $\alpha>2$ we have,
\begin{align*}
E(\exp(s\vert Z\vert))=\int_{-\infty}^{+\infty}\exp(s\vert z\vert)a_\alpha\exp(-\vert z\vert^\alpha)dz=2a_\alpha\mathcal{S},
\end{align*}
where
\begin{align*}
\mathcal{S}=\int_{0}^{+\infty}\exp(sz-z^\alpha)dz=\underbrace{\int_{0}^{1}\exp(sz-z^\alpha)dz}_{\mathcal{S}_1}+\underbrace{\int_{1}^{+\infty}\exp(sz-z^\alpha)dz}_{\mathcal{S}_2}.
\end{align*}
For $z\in[0,1]$ we have $\exp(-z^\alpha)\leq 1$ and so
\begin{align*}
\mathcal{S}_1\leq \int_{0}^{1}\exp(sz)dz=\frac{\exp(s)-1}{s}.
\end{align*}
For $z\geq 1$ we have $\exp(z^2-z^\alpha)<1$ and so
\begin{align*}
\mathcal{S}_2=\int_{1}^{+\infty}\exp(sz-z^2+z^2-z^\alpha)dz<\int_{1}^{+\infty}\exp(sz-z^2)dz<\sqrt{\pi}\exp\Big(\frac{s^2}{4}\Big),
\end{align*}
where the last inequality is obtained using Equation \eref{glap}. Finally, we obtain
\begin{align*}
\mathcal{S}<\frac{\exp(s)-1}{s}+\sqrt{\pi}\exp\Big(\frac{s^2}{4}\Big),
\end{align*}
and therefore 
\begin{align*}
E(\exp(s\vert Z\vert))<2a_\alpha\Big(\frac{\exp(s)-1}{s}+\sqrt{\pi}\exp\Big(\frac{s^2}{4}\Big)\Big).
\end{align*}

\hfill $\Box$ 

\subsection{Proof of Corollary \ref{lipschitz}}\label{prooflipschitz}
We suppose that the inequality (\ref{conshumed}) holds and we want to find an upper bound for $P\Big(\vert \phi(Z)-E(\phi(Z))\vert\geq u\Big)$. 
Using the Markov's inequality we have,
\begin{align}
\label{start}
P\Big(\vert \phi(Z)-E(\phi(Z))\vert >u\Big)&= P\Big(\exp(\lambda\vert \phi(Z)-E(\phi(Z))\vert ^2) > \exp(\lambda u^2)\Big),\nonumber\\
&\leq \exp(-\lambda u^2)E\Big(\exp(\lambda\vert \phi(Z)-E(\phi(Z))\vert ^2)\Big).
\end{align}

To demonstrate the result of the Theorem, it suffices to find an upper bound for the following quantity 
$$E\Big(\exp(\lambda\vert \phi(Z)-E(\phi(Z))\vert ^2)\Big).$$
Let $Z_1$ and $Z_2$ be two independent random variables distributed with the same law, then for all $u'>0$ we have:
\begin{equation}
\label{randmed}
P\Big(\vert Z_1-Z_2\vert > u'\Big)\leq P\Big(\vert Z_1-M(\phi(Z_1))\vert > \frac{u'}{2}\Big)+P\Big(\vert Z_2-M(\phi(Z_2))\vert > \frac{u'}{2}\Big).
\end{equation}
Furthermore, for all convex function $\psi$ we have:
\begin{align*}
E\Big(\psi(Z_1-E(Z_1))\Big) =\int \psi\Big(\int (z_1-z_2)dP(z_2)\Big)dP(z_1).
\end{align*}
Applying the Jensen's inequality we obtain then,
\begin{align}
\label{randmean}
E\Big(\psi(Z_1-E(Z_1))\Big)&\leq \int \Big(\int\psi (z_1-z_2)dP(z_2)\Big)dP(z_1)\nonumber,\\
&\leq E\Big(\psi(Z_1-Z_2)\Big).
\end{align}

Set $\psi(t)=\exp(\lambda t^2)$ for $\lambda>0$, $Z_1=\phi(Z)$ and $Z_2=\phi(Z')$. Since $\psi(t)$ is convex, then Equation (\ref{randmean}) gives:
\begin{align}
\label{randmeanconvex}
E\Big(\exp(\lambda\vert \phi(Z)-E(\phi(Z))\vert ^2)\Big)&\leq E\Big(\exp(\lambda(\phi(Z)-\phi(Z'))^2)\Big).
\end{align}
For all non-negative random variables $Z$ we have $E(Z)=\int_{[0,\infty)} P(Z\geq z)dz$. 
So, we obtain from Equation (\ref{randmeanconvex}):
\begin{align*}
E\Big(\exp(\lambda\vert \phi(Z)-E(\phi(Z))\vert ^2)\Big) \leq \int_0^\infty P\Big(\exp(\lambda(\phi(Z)-\phi(Z'))^2)>t\Big)dt.
\end{align*}
Using Equation (\ref{randmed}) and simple calculations leads to:
\begin{align}
E\Big(\exp(\lambda\vert \phi(Z)-E(\phi(Z))\vert ^2)\Big)&\leq \int_1^\infty P\Big(\vert\phi(Z)-\phi(Z')\vert >\sqrt{\frac{\log(t)}{\lambda}}\Big)dt,\nonumber\\
&\leq 2\int_1^\infty P\Big(\vert\phi(Z)-M(\phi(Z))\vert >\frac{1}{2}\sqrt{\frac{\log(t)}{\lambda}}\Big)dt.
\end{align}
In this step we can use the result in Equation (\ref{conshumed}), from which we obtain:
\begin{align}
\label{shuse}
E\Big(\exp(\lambda\vert \phi(Z)-E(\phi(Z))\vert ^2)\Big)\leq 2B\int_1^\infty \exp\Big(-\frac{\log(t)}{4\lambda A}\Big)dt,
\end{align}
and, therefore,
\begin{align}
\label{shuseres}
E\Big(\exp(\lambda\vert \phi(Z)-E(\phi(Z))\vert ^2)\Big)\leq \frac{8\lambda AB}{1-4\lambda A}\quad ,\quad \forall \lambda<\frac{1}{4A}.
\end{align}
The proof is complete by taking $\lambda =1/8A$.

\hfill $\Box$ 

\newpage
\bibliography{biblio}
\end{document}